\documentclass [12pt, a4paper, twoside]{article}

\usepackage {diagrams, a4, indentfirst, latexsym, theorem, QED, amssymb, epsfig}

\renewcommand{\qedsymbol}{\quad\Box}
\renewcommand{\TheWordProof}{\bf Proof.}
\newcommand{\doqed}{\pushright{$\qedsymbol$}}    
\newcommand\nothing[1]{\relax}

\theoremstyle{change}
\theorembodyfont{\it}
\newtheorem{theorem}{Theorem.}[subsection]
\newtheorem{proposition}[theorem]{Proposition.}
\newtheorem{lemma}[theorem]{Lemma.}
\newtheorem{corollary}[theorem]{Corollary.}
{\theorembodyfont{\upshape}
\newtheorem{definition}[theorem]{Definition.}
\newtheorem{remark}[theorem]{Remark.}

\newtheorem{construction}[theorem]{Construction.}

\newtheorem{observation}[theorem]{Observation.}
}

\newcommand{\hocolim}{\mathrm{hocolim}}
\newcommand{\Vis}{\mathit{Vis }}
\newcommand{\Invis}{\mathit{Inv }}
\newcommand{\Back}{\mathit{Back }}
\newcommand{\Front}{\mathit{Front }}
\newcommand{\Low}{\mathit{Low }}
\newcommand{\Up}{\mathit{Up }}
\newcommand{\pliso}{\iso_{PL}}
\newcommand{\bR}{{\mathbb R}}
\newcommand{\bZ}{{\mathbb Z}}
\newcommand{\bC}{{\mathbb C}}
\newcommand{\bN}{{\mathbb N}}
\newcommand{\bQ}{{\mathbb Q}}
\newcommand{\iso}{\cong}
\newcommand{\tensor}{\otimes}

\newcommand{\id}{\mathrm{id}}
\newcommand{\ie}{{\it i.e.}}
\newcommand{\hcoh}{\frak{hCoh}(P)}
\newcommand{\wchcoh}{\overline{\frak{hCoh}}(P)}
\newcommand{\hcohf}{\hcoh_{\mathrm{f}}}
\newcommand{\hcohhf}{\hcoh_{\mathrm{hf}}}
\newcommand{\wchcohhf}{\wchcoh_{\mathrm{hf}}}
\newcommand{\Knl}{K^{\mathrm{nl}}}
\newcommand{\Top}{\hbox{\bf Top}_*}
\newcommand{\dTop}{\hbox{-}\Top}
\newcommand{\C}{\hbox{\bf C}}
\newcommand{\Cf}{\C_{\mathrm{f}}}
\newcommand{\Chf}{\C_{\mathrm{hf}}}

\date {15.05.2007}

\begin{document}

\title{$K$-Theory of non-linear projective toric varieties}
\author{Thomas H\"uttemann}
\maketitle

\centerline {\it Queen's University Belfast, Department of Pure
Mathematics}

\centerline {\it Belfast BT7~1NN, Northern Ireland, UK}

\centerline {e-mail: \texttt{t.huettemann@qub.ac.uk}}

\vglue 2\bigskipamount \hrule \medskip

{\footnotesize We define a category of
  quasi-coherent sheaves of topological spaces on
  projective toric varieties and prove a splitting result for its algebraic
  $K$-theory, generalising earlier results for projective spaces. The
  splitting is expressed in terms of
  the number of interior lattice points of dilations of a polytope
  associated to the variety. The proof uses combinatorial and geometrical results
  on polytopal complexes. The same methods also give an
  elementary explicit calculation of the cohomology groups of a projective
  toric variety over any commutative ring.
\smallskip


}
\medskip \hrule \vglue 2\bigskipamount

Let $X$ denote a scheme with a preferred covering by open affine subschemes.
A vector bundle on~$X$ can be described by a collection of finitely generated
projective modules, one for each open affine of the chosen covering, and
``restriction'' maps between them, satisfying a certain gluing condition. For
toric varieties it is possible to define analogous topological objects,
replacing rings by monoids, modules by topological spaces, and weakening the
gluing condition to a homotopy invariant condition. This program has been
carried out by the author for projective spaces in~\cite{H-Proj} where it was
shown that the algebraic $K$-theory of the resulting category of ``non-linear
sheaves'' splits into $n+1$ copies of \textsc{Waldhausen}'s $K$-theory space
$A(*)$. The aim of the present paper is to generalise this splitting result to
arbitrary projective toric varieties, thereby revealing much of the
combinatorial content of the earlier result explicitly.

This paper can also be understood as an attempt to describe toric varieties
over ``brave new rings'', replacing (commutative) rings by ring spectra. The
combinatorial structure of toric varieties is rigid enough to allow a
treatment with techniques from unstable homotopy theory (using spaces, not
spectra). It is not clear what a toric variety should be in that context, but
we can nevertheless define quasi-coherent sheaves on such a variety, called
non-linear sheaves. This category carries enough structure to define, for
example, algebraic $K$-theory, just as the $K$-theory of a ring~$R$ can be
defined in terms of a category of $R$-modules.

\medbreak

An $n$-dimensional polytope with integral vertices defines a
projective toric variety~$X_P$ (its construction is reviewed in
\S\ref{subsec:barr-cones-proj}), equipped with an ample
(equivariant) line bundle~$\mathcal{O}(1)$. (It can be shown that
any projective toric variety over~$\bC$ equipped with an ample
(equivariant) line bundle arises in this way.)  We denote the
algebraic $K$-theory of non-linear sheaves on~$X_P$ by~$\Knl(X_P)$.
The following is the main theorem of the paper:

\bigbreak \noindent
{\bf Theorem~\ref{thm:splitting}.} {\it Let $P
\subset \bR^n$ be a
  polytope with integral vertices, and assume that $P$ has non-empty
  interior. Let $k$ denote the number of integral
  roots of its \textsc{Ehrhart\/} polynomial (cf.~Theorem~\ref{thm:Ehrhart}); possibly $k=0$. Then there is
  a homotopy equivalence\goodbreak
  \[\underbrace{A(*) \times \ldots \times A(*)}_{(k+1) \hbox{ factors}}
     \times \Knl (X_P)^{[k]}
     \rTo^\sim \Knl (X_P)\]
  where $\Knl (X_P)^{[k]}$ is the algebraic $K$-theory of the
  category of those non-linear sheaves~$Y$ on~$X_P$ which have $\Gamma (Y(i))
  \simeq *$ for all $0 \leq i \leq k$.}
\medbreak

\noindent In fact, $k\geq 0$ is minimal among integers $j \geq 0$
such that the dilated polytope $(j+1)P$ has a lattice point in its
interior. The map in the theorem is induced by the assignment
$$ (K_0,\,K_1,\,\ldots,\,K_k) \mapsto \bigvee_{i=0}^k K_i \wedge \mathcal{O}_P
(-i) \ .$$ Here the $K_i$ are pointed topological spaces,
$\mathcal{O}_P$ is the non-linear analogue of the structure sheaf
on~$X_P$, and $\mathcal{O}_P(i)$ is its $i$th twist
(\S\ref{subsec:non-linear-sheaves}). The functor $\Gamma$ is the
total cofibre functor (\S\ref{subsec:total-cofibres}), a substitute
for the global sections functor and its derived functors in
algebraic geometry.

A similar splitting result should hold for the algebraic $K$-theory
of projective toric varieties over~$\bC$. To explain the passage to
the ``linear'' world, note that by taking free $\bC$-vector spaces
the non-linear sheaves $\mathcal{O}_P (j)$ give rise to the usual
twisting sheaves $\mathcal{O} (j) = \mathcal{O}(1)^{\otimes j}$ on
the $\bC$-scheme $X_P$. The meaning of the number~$k$ in the theorem
is that $H^i (X_P, \mathcal{O} (-k)) = 0$ for all $i \geq 0$, but
$H^n (X_P, \mathcal{O}(-k-1)) \not= 0$. It turns out that the
obstruction to a further splitting of $K(X_P)$ is the non-vanishing
of the cohomology of $\mathcal{O}(-k-1)$, or equivalently, the
presence of a lattice point in the interior of $(k+1)P$.

The total cofibres of the non-linear sheaves $\mathcal{O}_P (j)$ exhibit the
same behaviour as their linear counterparts $\mathcal{O}(j)$ as is
demonstrated by the explicit calculations in \S\ref{subsec:total-cofib-can}.
In particular, the obstruction for splitting off a further copy of $A(*)$ in
Theorem~\ref{thm:splitting} is the non-triviality of the total cofibre of
$\mathcal{O}_P (-k-1)$. The similarity between total cofibres and sheaf
cohomology is not coincidental as is explained in \cite{H-Cofibres}.

If $P$~is an $n$-dimensional standard simplex (\ie, if $P$ is a lattice
simplex with volume $1/n!$), then $k=n$, and the variety $X_P$ is the
$n$-dimensional projective space equipped with the usual twisting sheaf
$\mathcal{O}(1)$. In this case, it can be shown that $\Knl (X_P)^{[n]} \simeq
*$, so the above theorem reduces to the known splitting of~\cite{H-Proj}.

The proof of the splitting result relies on explicit computations of certain
homotopy colimits of ``geometrically defined'' diagrams.  We review polytopal
complexes in \S\ref{sec:excision}, and prove a result relating the nerve of a
subset of a polytopal complex to its underlying space. Results of this type, often
well-known for simplicial complexes,
are part of the general toolkit for dealing with homotopy types of nerves of
posets \cite{WZZ-hocolims}; the given version is slightly
more general than needed for this paper. Other geometrical issues are
discussed in \S\ref{subsec:BdComplex}: Links, stars and visibility
subcomplexes of the boundary complex of a polytope are introduced and
examined in detail.  Although quite elementary, it will be important for the
sequel to provide explicit descriptions throughout.

In \S\ref{sec:non-linear-sheaves} we review the description of quasi-coherent
sheaves on projective toric varieties by diagrams of modules. By analogy, a
notion of ``non-linear sheaves'' is introduced
(\S\ref{subsec:non-linear-sheaves}). Twisting sheaves and tensor products are
also defined by analogy. To obtain a homotopically meaningful analogy of
global sections, we use the ``total cofibre'' construction
of~\S\ref{subsec:total-cofibres}. The vanishing criterion for total cofibres
makes use of the combinatorial results from~\S\ref{sec:polytopal-complexes}.
Next, we calculate total cofibres of twisting sheaves
(\S\ref{subsec:total-cofib-can}); here the material from
\S\ref{sec:polytopal-complexes} is used heavily again.  We show that the same
techniques also lead to an elementary computation of the cohomology groups
$H^r (X; \mathcal{O} (k))$ over any commutative ring~$R$ where $X = X_P$ is
the toric variety defined by~$P$.  (Standard references for toric geometry
seem to miss an explicit combinatorial treatment of negative twists. Note also
that the given treatment does not use \textsc{Serre\/} duality to deal with the
case of non-ample line bundles but yields a direct identification of a basis
of the unique non-trivial cohomology module. See
Remark~\ref{rem:compute_cohomology} for pointers to an algebro-geometric
approach.)  Comparison with \S\ref{subsec:total-cofib-can} shows that the
total cofibre construction captures not only global sections, but higher
cohomology groups as well.

Finally, \S\ref{sec:algebraic-k-theory} is concerned with $K$-theoretical
issues. Following a brief discussion of finiteness conditions for non-linear
sheaves (which are also the subject of the paper \cite{H-Finiteness}) we
define their algebraic $K$-theory and prove the splitting result.

\vfill

\section{On Polytopal Complexes}
\label{sec:polytopal-complexes}

\subsection{Complexes and order filters}
\label{sec:excision}

A polytope~$P$ is the convex hull of a finite set of points in~$\bR^n$. We
write $F \leq P$ if~$F$ is a face of~$P$; this includes the case of improper
faces $F = P$ and $F=\emptyset$.

\begin{definition}
A non-empty finite collection~$K$ of non-empty polytopes in
some~$\bR^n$ is called a {\it polytopal complex\/} if the following
conditions are satisfied:
\begin{enumerate}
        \item If $ F \in K $ and $ \emptyset \not= G \leq F $, then $ G \in K $
        \item For all $ F, G \in K $, the intersection $ F \cap G $
        is a (possibly empty) face of~$F$ and~$G$.
\end{enumerate}
A subset $ L \subseteq K $ of a polytopal complex is called
an {\it order filter\/} if for all $ F \in L $ and $ G \in K $
with $ F \leq G $, we have $ G \in L $.
A subset $ L \subseteq K $ of a polytopal complex is called
a {\it subcomplex\/} of~$K$ if $L$~is a polytopal complex.
\end{definition}

Important examples of polytopal complexes are the complex~$F(P)_0$ of non-empty
faces of a polytope~$P$, and its subcomplex~$F(P)^1_0$ of non-empty
proper faces of~$P$ (sometimes called boundary complex of~$P$).

The intersection of two subcomplexes, if non-empty, is a subcomplex. The
(set-theoretic) complement of a subcomplex is an order filter.

\begin{definition}
Suppose $K$~is a polytopal complex, and $L$~is a non-empty subset
of~$K$. We call $ |L| := \bigcup_{F \in L} F $ the {\it
realisation\/} or the {\it underlying space\/} of~$L$. The {\it
combinatorial closure\/} of~$L$ in~$K$ is the set of all polytopes
in~$L$ and their non-empty faces:
        $$ \overline L := \left\{ F \in K \mid \exists G \in L \colon F \leq G
        \right\} \ . $$
\end{definition}

The combinatorial closure of~$L$ in~$K$ is a complex, and we have $ |L| =
|\overline L|$.  If~$P$~is an $n$-dimensional polytope, we have
$PL$-homeomorphisms $ |F(P)_0| = P \pliso B^n $ and $ |F(P)_0^1| = \partial P
\pliso S^{n-1} $.

\bigbreak

Note that a complex~$K$ is naturally a partially ordered set with order given
by inclusion of faces. Hence we can view any non-empty subset $L \subseteq K$
as a category with morphisms corresponding to inclusion of faces. Its nerve
$NL$~is an abstract simplicial complex; a $k$-simplex is a strictly increasing
sequence \hbox {$ [F_0 < F_1 < \ldots < F_k] $} of polytopes in~$L$. For each polytope
$ F \in L$ there is a corresponding vertex~$[F]$ of~$NL$. We denote the
geometric realisation of~$NL$ by~$|NL|$; this space is called the {\em
classifying space\/} of~$L$.

For $F \in K$ let~$\hat F$ denote its barycentre. Define a map $ \alpha \colon
|NL| \rTo |K| $ by sending the zero-simplex $[F] \in NL$ to the point $ \hat F
\in |K|$ and extending linearly over simplices. This map is an embedding and
thus allows us to view the abstract simplicial complex~$NL$ as a simplicial
complex, \ie, a polytopal complex consisting of simplices.

\begin{lemma}
\label{Triangulation}
Suppose $K$~is a polytopal complex. The simplicial complex~$NK$ is the
barycentric subdivision of~$K$. The map $ \alpha \colon |NK| \rTo |K| $ is
a $PL$-homeo\-morphism, and the pair $(NK,\ \alpha)$ is a triangulation of~$|K|$
(\cite{RS-PL}, p.~17).
\qed
\end{lemma}

\begin{definition}
Let~$K$ denote a polytopal complex, and fix $A \in K$.
\begin{enumerate}
        \item The {\it (open) star\/} of~$A$ is $ st (A) := \{ F \in K | A \leq F \}$.
        \item The {\it closed star\/} of~$A$ is defined as the combinatorial closure
        of $st (A)$. Explicitly, $ \overline {st} (A) = \{ F \in K | \exists G \in st(A):
        F \leq G \}$.
        \item The {\it (closed) antistar\/} of~$A$ is
        $ \overline {ast} (A) := K \setminus st (A) $.
        \item The {\it open antistar\/} of~$A$ is $ ast (A) := K \setminus \overline {st} (A) $.
        \item The {\it link\/} of~$A$ is $ lk (A):= \overline {st} (A) \cap \overline {ast} (A) $.
\end{enumerate}
If we have to emphasise the complex~$K$ we write $lk_K (A)$ instead of
$lk(A)$, and similar for the other expressions.
\end{definition}

The sets just defined are combinatorial rather than geometric in nature; for
example, $|st(A)| = |\overline {st} (A)|$, but as sets, $st(A)$ and $\overline
{st} (A)$ usually differ.

The open star and open antistar are order filters. The closed antistar, if
non-empty, is a subcomplex of~$K$.  The link $lk (A)$ can also be described by
$lk(A) = \overline{st}(A) \setminus st(A) = \overline{ast}(A) \setminus
ast(A)$. Geometrically the link of~$A$ consists of those polytopes that are
``visible'' inside~$|K|$ from the barycentre~$\hat A$ of~$A$ but do not
contain~$\hat A$.

Note that the open antistar may be empty while the closed antistar is not.
Thus in general $\overline{ast}(A)$~is not the combinatorial closure
of~$ast(A)$, but see below for the case of manifolds. (For example, consider
the complex $ K = F(P)_0 $. Then $\overline {ast} (P) = F(P)_0^1 = lk (P)$
while $ast(P) = \emptyset$.)

Let $L$ denote the stellar subdivision of~$K$ at~$\hat A$
\cite[Definition~III.2.1]{Ewald-CCAG}. By definition, $L$~is obtained from~$K$
by removing $st_K (A)$ and adding the cones from $\hat A$ on polytopes in~$lk_K(A)$.
Then it is easy to see that $ |lk_K (A)| = |lk_{L} (\hat A)|
$ and $ |\overline {st}_K (A)| = |\overline {st}_{L} (\hat A)|$, and this
agrees with the definition of link and star in~\cite{RS-PL}, p.~20. Thus we
see that $|st(A)|$~is a (topological) neighbourhood of~$\hat A$ in~$|K|$; it
is the cone from~$\hat A$ on~$|lk(A)|$.  If $K$~is an $m$-dimensional $PL$
manifold (possibly with boundary), $|st(A)| \pliso B^m$, while $|lk(A)|$ is
$PL$ homeomorphic to $S^{m-1}$ if $\hat A \in \mathrm {int}\, |K|$ and to
$B^{m-1}$ if $\hat A \in \partial |K|$ by \cite[Exercise~2.21~(1)]{RS-PL}.  Moreover, if
$\hat A \in \mathrm {int}\, |K|$ we know that $|lk (A)| = \partial |st(A)|$ is
the boundary sphere of the ball~$|st(A)|$.

\begin{lemma}
\label{AstClosure}
  Suppose~$|K|$ is an $m$-dimensional $PL$~manifold without boundary. Then $ast(A) \not=
  \emptyset$, and the closed antistar of~$A$ is a combinatorial
  closure of~$ast(A)$. Moreover, $|ast(A)|$ is the closure of the complement
  $|K| \setminus |st(A)|$, and $|lk (A)| \pliso S^{m-1}$ is the boundary of
  both $|st(A)|$ and~$|ast(A)|$.
\end{lemma}

\begin{proof}%
Note first that $ ast (A) = \emptyset $ implies $\overline {st} (A) = K$. But
$|\overline {st} (A)|$ is a ball since~$|K|$ is a $PL$ manifold \cite[Exercise~2.21~(1)]{RS-PL}.
This contradicts the assumption that $K$~has no boundary. Thus
necessarily $ast (A) \not= \emptyset$.

Since $\overline {ast} (A)$ is a complex, the combinatorial closure of~$ast
(A)$ is contained in the closed antistar. Conversely,
given an element $F \in \overline {ast} (A) \setminus ast (A)$ we have to show
that $F$~is the face of some $G \in ast (A)$. Suppose such~$G$ does not
exist. Then $st (F) \subseteq \overline {st} (A)$. Moreover, $|st(F)|$~is known to be
a ball with~$\hat F$ in its interior, which shows that $\hat F$~is an interior
point of the ball $|\overline {st} (A)|$. But since $F \in lk(A)$ we know $\hat F
\in \partial |\overline {st} (A)|$, a contradiction.---%
The other assertions are obvious.
\end{proof}

\medbreak

A similar argument shows more generally:

\begin{lemma}
\label{Anticomplex}
  Let $L$ be a subcomplex of~$K$. Suppose $|K|$ is an $m$-dimensional $PL$
  manifold without boundary, and $|L|$ is an $m$-dimensional $PL$
  manifold with boundary. Then
  $C:= K \setminus L$ is non-empty, and $|C|$~is the closure of the complement
  $|K| \setminus |L|$. Let~$B$ denote the intersection of~$L$ and~$\overline
  C$. Then $|B|$~is the boundary of both~$|L|$ and~$|C|$.
  \qed
\end{lemma}

The following lemma shows how to connect data useful for homotopy theory
(the nerve of a certain category) with geometric data (subspaces of a $PL$
manifold). Similar results are well-documented for simplicial complexes;
extension to polytopal complexes can be achieved by passage to barycentric
subdivisions. Our proof utilises the ``simplicial neighbourhood theorem'' of
\textsc{Rourke\/} and \textsc{Sanderson\/} \hbox{\cite[Theorem~3.11]{RS-PL}}.

\begin{lemma}
  Let~$K$ denote a polytopal complex, and let $C \subseteq K$ be an order
  filter.  Assume $|K|$ is an $m$-dimensional $PL$ manifold without boundary
  and $|K \setminus C|$ is an \hbox {$m$-dimensional} $PL$ manifold with
  boundary.  Then $|C|$ is a regular neighbourhood of~$|NC|$ in~$|K|$, hence $|NC|$ is
  a deformation retract of~$|C|$.
\end{lemma}

\begin{proof}%
  The map~$\alpha$ from~\ref{Triangulation} allows us to consider~$|NC|$ as a
  subspace of~$|C|$. Given that~$|C|$ is a regular neighbourhood of~$|NC|$,
  the collapsing criterion \cite[Corollary~3.30]{RS-PL} shows that $|C|$
  collapses to $|NC|$, thereby proving the proposition.

Define $B:= \overline C \setminus C$.  Using Lemma~\ref{Triangulation} we see
that $|C| = |\overline C| \pliso |N\overline C|$. Moreover, $N C$ and $N B$ are
simplicial subcomplexes of $N \overline C$.  By the ``simplicial neighbourhood
theorem'' \cite[Theorem~3.11]{RS-PL}, it thus suffices to prove the
following assertions:
\begin{enumerate}
\item \label{bd} $|N\overline C|$ is a compact $PL$ manifold with
  boundary, and $N B$ is a triangulation of $\partial |N\overline C|$
\item \label{full} $N C$ is a full subcomplex of~$N \overline C$
\item \label{comp} $N B$ is the simplicial complement of $N C$
inside $N \overline C$
\item \label{nbhd} $N \overline C$ is the simplicial neighbourhood of
$N C$ in $N \overline C$
\item \label{int} $|NC|$ lies in the interior of $|N\overline C|$
\end{enumerate}

(\ref{bd}): Lemma~\ref{Anticomplex} implies that $|C|$~is a $PL$ manifold with
boundary~$|B|$ (compactness is automatic since all our complexes are finite).
By~\ref{Triangulation} there are homeomorphisms $|C| \pliso |N\overline C|$ and
$|B| \pliso |NB|$. Thus $NB$ is a triangulation of $\partial |C| = |B|$.

(\ref{full}): A $k$-simplex~$F$ of $N \overline C$ is a chain of polytopes
$$ F = [A_0 < A_1 < \ldots < A_k] $$
with $A_i \in \overline C$. Assume the boundary of~$F$
is contained in $ N C $. Then in particular all its vertices $[A_i]$ are in $NC$,
\ie, $A_i \in C$, hence $F$~is an element of~$N C$ by definition of the nerve.
By \cite{RS-PL}, Exercise~3.2 this implies assertion~\ref{full}.

(\ref{comp}): The simplicial complement of $N C$ inside $N \overline C$ is, by
definition, the set
$$ \{ F \in N \overline C \mid F \cap |NC| = \emptyset \}\ .$$
Let~$F = [A_0 < \ldots < A_k]$ be a $k$-simplex in $N \overline C$. Then $F \cap |NC|
\not= \emptyset$ if and only if there is a simplex $G = [B_0 < \ldots <
B_\ell] \in N C$ with $F \cap G \not= \emptyset$. But $F \cap G$ is also
a face of~$F$ and~$G$. In particular, $F$ and~$G$ have non-empty intersection
if and only if they have a common vertex $A_i = B_j \in C$. This shows
that $ F \cap |NC| = \emptyset $ if and only if no~$A_i$ is in~$C$,
\ie, if and only if~$F \in N B = N(\overline C \setminus C)$.

(\ref{nbhd}): The simplicial neighbourhood of $N C$ in $N \overline
C$ is, by definition, the set
$$ T:= \{ F \in N\overline C \mid \exists G \in N\overline C:\ F \leq G
\ \mathrm {and}\ G \cap |NC| \not= \emptyset \} \ .$$
Let $ F = [A_0 < \ldots < A_k] $ be a $k$-simplex of $N \overline
C$. If $A_k \in C$ we have $F \cap |NC| \not= \emptyset$ by the arguments
in~(\ref{comp}), hence $ F \in T $. Otherwise there exists $ A_{k+1} \in C $
with $ A_k < A_{k+1} $ by definition of the combinatorial closure.  Then
$F$~is a face of $G := [A_0 < \ldots <A_k < A_{k+1}] $ and $ G \cap |NC| \not=
\emptyset $ by construction, thus $ F \in T $.

(\ref{int}): We show that $|NC| \cap \partial |N\overline C| =
\emptyset $. Recall that $\partial |N\overline C| = |B|$. Let~$F
\in N \overline C$ be given. By arguments similar to those
of~(\ref{comp}), applied to the complex~$N B$, we see that $ F \cap |B| $
is non-empty if and only if the chain $ A_0 < \ldots < A_k $
representing~$F$ satisfies $A_0 \in B$. But that cannot happen for $F \in
N C$.
\end{proof}

\begin{corollary}
\label{RegularNbhdSpheres}
  Let $K$ denote a polytopal complex, and let $C \subseteq K$ be an order filter.
  Suppose $|K|$ is an $m$-dimensional $PL$ sphere, and $|K \setminus C|$ is an
  \hbox{$m$-dimensional} ball. Then $|NC|$ is contractible.
\end{corollary}

\begin{proof}%
  From Lemma~\ref{Anticomplex} and \cite[Corollary~$3.13_m$]{RS-PL}, we know
  that~$|C|$ is a $PL$ ball, hence is contractible. Consequently its
  deformation retract~$|NC|$ is contractible as well.
\end{proof}


\subsection{The boundary complex of a polytope}
\label{subsec:BdComplex}

We restrict attention to the special case of the boundary complex of an
\hbox{$n$-dimensional} polytope $P \subseteq \bR^n$. Its realisation is
$\partial P$, thus it is a $PL$ sphere of dimension \hbox{$n-1$}. In order to
apply Corollary~\ref{RegularNbhdSpheres} we need to construct
\hbox{``interesting''} $(n-1)$-balls inside~$\partial P$. One class of
examples is given by the closed stars which can be characterised by purely
combinatorial means. We also discuss examples given by subsets of faces
satisfying certain geometric conditions.

\subsubsection*{Links, stars and antistars}

For a polytope~$P$, the set $F(P)_0^1$ of non-empty proper faces of~$P$
is a polytopal complex, called {\it boundary complex\/} of~$P$.  Links and
antistars admit convenient combinatorial descriptions in this case.

The set~$F(P)$ of all faces of~$P$ (including~$P$ and~$\emptyset$) is known to
be a finite graded lattice \cite[Theorem~2.7]{Ziegler-Polytopes}. We write $F
\vee G$ for the join of~$F$ and~$G$ in~$F(P)$; it is the smallest face of~$P$
containing $F \cup G$. Links, stars and antistars are computed in the complex
$K = F(P)_0^1$ unless indicated otherwise; in particular, the star of a proper
face of~$P$ will not contain~$P$ itself.

\begin{lemma}
\label{CombDesc}
{\em (Combinatorial description of star, link and antistar.)}\\
Let~$A$ denote a proper non-empty face of~$P$.
\begin{enumerate}
\item $ \overline{st}(A) = \{ F \in F(P)_0^1 | F \vee A \not= P\} $
\item $ lk(A) = \{ F \in F(P)_0^1 | F \vee A \not= P \ \&\ A \not\leq F \}$
\item $ ast(A) = \{ F \in F(P)_0^1 | F \vee A = P \} $
\end{enumerate}

\end{lemma}

\begin{proof}%
To prove~(1), suppose $F\in F(P)_0^1$ satisfies $ F \vee A \not= P$. Then $F
\vee A \in F(P)_0^1 $, and from $ F \leq F \vee A \geq A $
we get $ F \in \overline{st} (A) $.
Conversely, if~$F$ is an element of the closed star of~$A$, we find a proper
face~$G$ of~$P$ with $ F \leq G \geq A $. But then $ F \vee A \leq G \vee A = G \not= P $.

Assertions~(2) and~(3) follow immediately from~(1).
\end{proof}

\begin{corollary}
\label{CombStarLinkCor}
Let~$A$ be a proper non-empty face of~$P$.
If $ B \in st(A) \setminus \{A\} $, we have $ A \in lk(B) $ and
$$ \overline{st}_{lk(B)}(A) = \{ F \in lk(B) | B \not\subseteq F \vee A \} $$
where $\overline{st}_{lk(B)}(A)$ denotes the closed star of~$A$ in the
polytopal complex~$lk(B)$.
\end{corollary}

\begin{proof}%
By hypothesis $B \supset A$, hence $A \in lk(B)$.
Suppose we have an element $ F \in \overline{st}_{lk(B)}(A) $.
By definition of the closed star, there is a $ G \in lk(B) $
with $ F \leq G \geq A $. But then
$$ F \vee A \subseteq G \vee A = G \ \in lk(B) \ .$$
Since $ B \not\subseteq G $ by definition of the link, this implies
$ B \not\subseteq F \vee A $.

Conversely, given $ F \in lk(B) $ with $ B \not\subseteq F \vee A $,
we know that $ F \vee A \subseteq F \vee A \vee B = F \vee B \not= P $, thus $ F \vee A \in lk(B) $.
From $ F \leq F \vee A \geq A $ we conclude $ F \in \overline{st}_{lk(B)}(A)$.
\end{proof}

\subsubsection*{Visible and invisible faces}
\label{SecVisInvis}

\begin{definition}
  A face $ F \in F(P)_0^1 $ is called {\it visible\/} from the point $x \in
  \bR^n \setminus P$ if $ [p,x] \cap P = \{p\} $ for all $ p \in F $. (Here
  $[p,x]$ denotes the line segment between~$p$ and~$x$.) Equivalently, $F$~is
  visible if $p + \lambda (x-p) \notin P$ for all points $p \in F$ and real
  numbers $\lambda > 0$.  We denote the set of visible faces by~$\Vis(x)$; its
  complement $ \Invis(x) := F(P)_0^1 \setminus \Vis(x) $ is the set of {\it
    invisible\/} faces. Let $ \overline {\Invis(x)} $ denote the combinatorial
  closure of~$\Invis(x)$, and define $ \partial \Invis(x) := \overline
  {\Invis(x)} \setminus \Invis(x) = \overline {\Invis(x)} \cap \Vis(x) $.
\end{definition}

\begin{figure}[ht]
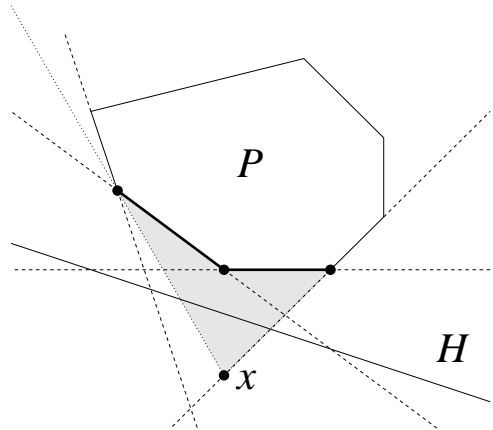

  \begin{center}
    \input visible.pstex_t
    \caption{Visible faces}
    \label{fig:visible}
  \end{center}
\end{figure}

\begin{lemma}
\label{VisibleFacets}
  A facet~$F$ of~$P$ is visible from~$x$ if and only if
  $x$~and $\mathrm {int}\, P$ are on different sides of the affine hyperplane
  spanned by~$F$. A proper non-empty face of~$P$ is visible if and only if
  it is contained in a visible facet of~$P$.
\qed
\end{lemma}

In particular, the sets~$\Vis(x)$ and~$\Invis(x)$ are non-empty.  Since a face of a
visible face is visible itself, $\Vis(x)$ and~$\partial \Invis(x)$ are subcomplexes
while $\Invis(x)$~is an order filter. If~$x$ is beyond~$F$ in the sense
of~\cite{Ziegler-Polytopes}, p.~78, the set of visible faces~$\Vis(x)$ coincides
with the closed star of~$F$.

\begin{proposition}
\label{VandIareBalls}
  \begin{enumerate}
  \item \label{VisBall} There is a $PL$ homeomorphism $|\Vis(x)| \iso B^{n-1}$.
  \item \label{IisBall} There is a $PL$ homeomorphism $|\Invis(x)| \iso B^{n-1}$.
  \item $ |\Vis(x)| \cap |\Invis(x)| = |\partial \Invis(x)| $ is the
    common boundary of both~$|\Vis(x)|$ and~$|\Invis(x)|$, hence is
    $PL$-homeomorphic to~$S^{n-2}$.
  \end{enumerate}
\end{proposition}

\begin{proof}%
Applying a translation if necessary we may assume~$x = 0$. For statement~(1), let~$H$
be any hyperplane separating~$0$ and~$P$ (Fig.~\ref{fig:visible}). Let~$C$ denote the cone (with
apex~$0$) on~$P$. Then $C$~is a pointed polyhedral cone, hence $C \cap H$ is a
$PL$ ball \cite[Theorem~V.1.1]{Ewald-CCAG}. Projection along~$C$ provides a
homeomorphism $|\Vis(x)| \iso C \cap H$. By the ``pseudo radial projection''
technique (\cite{RS-PL}, proof of Lemma~2.19) this can be modified to
give a $PL$ homeomorphism.

Statements~(2) and~(3) follow from Lemma~\ref{Anticomplex} and the fact that
the closure of the complement of a (full dimensional) $PL$ ball inside a $PL$
sphere is a $PL$ ball itself (\cite{RS-PL}, Corollary~$3.13_n$).
\end{proof}

\bigbreak

\begin{corollary}
\label{BIcontractible}
  The classifying space of $\Invis(x)$ is a deformation retract of $|\Invis(x)|$. In
  particular, $|N\Invis(x)|$ is contractible.
\end{corollary}

\begin{proof}%
This follows from Corollary~\ref{RegularNbhdSpheres} applied to $K = F(P)_0^1$ and
$C = \Invis(x)$, using Proposition~\ref{VandIareBalls}~(1).
\end{proof}

\subsubsection*{Front and back faces}

\begin{definition}
\label{DefFrontBack}
  A face $F \in F(P)_0^1$ is called a {\it back face\/} with respect to the
  point $x \in \bR^n \setminus \mathrm {int}\,P$
  if for all points $p \in F$ and all real numbers $\lambda > 0$ we have $p +
  \lambda(p-x) \notin P$. The set of back faces is denoted by~$\Back(x)$; its
  complement $\Front(x):= F(P)_0^1 \setminus \Back(x)$ is the set of {\it
    front faces}. Let $ \overline {\Front(x)} $ denote the combinatorial
  closure of~$\Front(x)$, and define $\partial \Front(x):= \Back(x) \cap
  \overline{\Front(x)} $.
\end{definition}

\begin{figure}[h]
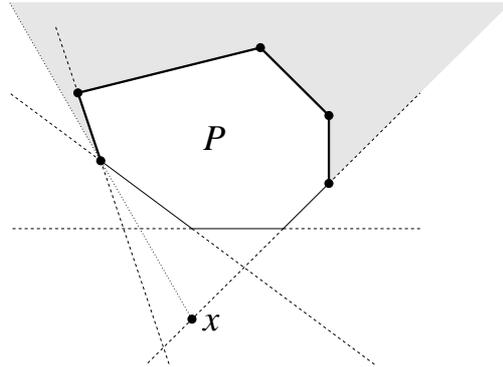

  \begin{center}
    \input back.pstex_t
    \caption{Back faces}
    \label{fig:back}
  \end{center}
\end{figure}

\begin{lemma}
\label{LemmaFrontBack}
  Suppose~$F$ is a facet of~$P$. Then~$F$ is a back face with respect to~$x$
  if and only if~$x$ and~$\mathrm {int}\, P$ are on the same side of the
  affine hyperplane spanned by~$F$. A proper non-empty face~$F$ of~$P$ is a
  back face if and only if it is contained in a facet of~$P$ which is a back face.
\qed
\end{lemma}

In particular, the sets~$\Back(x)$ and~$\Front(x)$ are non-empty.  Since a face of a
back face is a back face itself, $\Back(x)$ and~$\partial \Front(x)$ are subcomplexes
while $\Front(x)$~is an order filter.

By arguments similar to the ones used for the case of visible faces, we can show:

\begin{proposition}
\label{FandBareBalls}
  \begin{enumerate}
  \item There is a $PL$ homeomorphism $|\Back(x)| \iso B^{n-1}$.
  \item There is a $PL$ homeomorphism $|\Front(x)| \iso B^{n-1}$.
  \item $ |\Front(x)| \cap |\Back(x)| = |\partial \Front(x)| $ is the
    boundary of both~$|\Front(x)|$ and~$|\Back(x)|$, hence is $PL$-homeomorphic
    to~$S^{n-2}$.
    \qed
  \end{enumerate}
\end{proposition}

\begin{corollary}
\label{BFcontractible}
  The classifying space of $\Front(x)$ is a deformation retract of $|\Front(x)|$. In
  particular, $|N\Front(x)|$ is contractible.
\qed
\end{corollary}

\subsubsection*{Upper and lower faces}

\begin{definition}
\label{LowUpDef}
  A face $F \in F(P)_0^1$ is called a {\it lower face\/} with respect to the
  direction $x \in \bR^n \setminus \{0\}$
  if for all points $p \in F$ and all real numbers $\lambda > 0$ we have $p -
  \lambda x \notin P$. The set of lower faces is denoted by~$\Low(x)$; its
  complement $\Up(x):= F(P)_0^1 \setminus \Low(x)$ is the set of {\it upper
    faces}. Let $ \overline {\Up(x)} $ denote the combinatorial closure
  of~$\Up(x)$, and define $\partial \Up(x):= \Low(x) \cap \overline{\Up(x)} $.
\end{definition}

\begin{figure}[h]
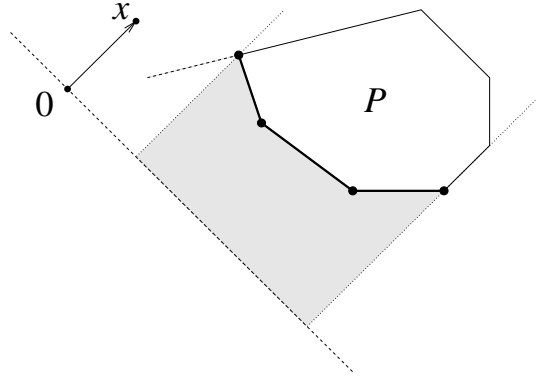

  \begin{center}
    \input lower.pstex_t
    \caption{Lower faces}
    \label{fig:lower}
  \end{center}
\end{figure}

\begin{lemma}
  Suppose~$F$ is a facet of~$P$ with inward pointing normal
  vector~$v$. Then~$F$ is a lower face with respect to~$x$ if and only if
  $\langle x, \, v \rangle > 0$. A proper non-empty face of~$P$ is a lower
  face if and only if it is contained in a facet of~$P$ which is a lower face.
\qed
\end{lemma}

In particular, the sets~$\Low(x)$ and~$\Up(x)$ are non-empty.  Since a face of a
lower face is a lower face itself, $\Low(x)$ and~$\partial \Up(x)$ are subcomplexes
while $\Up(x)$~is an order filter.

By arguments similar to the ones used for the case of visible faces, we can show:

\begin{proposition}
\label{UandLareBalls}
  \begin{enumerate}
  \item There is a $PL$ homeomorphism $|\Low(x)| \iso B^{n-1}$.
  \item There is a $PL$ homeomorphism $|\Up(x)| \iso B^{n-1}$.
  \item $ |\Low(x)| \cap |\Up(x)| = |\partial \Up(x)| $ is the common boundary of
    both~$|\Low(x)|$ and~$|\Up(x)|$, hence is $PL$-homeomorphic to~$S^{n-2}$.
    \qed
  \end{enumerate}
\end{proposition}

\begin{corollary}
\label{BUcontractible}
  The classifying space of $\Up(x)$ is a deformation retract of $|\Up(x)|$. In
  particular, $|N\Up(x)|$ is contractible.
\qed
\end{corollary}


\section{Non-Linear Sheaves and Total Cofibres}
\label{sec:non-linear-sheaves}

\subsection{Equivariant spaces}
\label{subsec:EquivSpaces}

Before describing quasi-coherent sheaves on projective toric varieties
we introduce some terminology concerning topological spaces.  Let $M$
denote an abelian pointed monoid (\ie, we have elements $*,0 \in M$
such that $0+m=m$ and $*+m=*$ for all $m \in M$).  Any abelian monoid
can be made into a pointed monoid by artificially adding a disjoint
basepoint~$*$. We consider $M$ as a discrete topological space with
$*$ as base point. The category of pointed topological spaces with a
right (base point preserving) action of~$M$ will be denoted $M\dTop$.
The $M$-equivariant $n$-cell is the space $D^n_+ \wedge M$, its
boundary is $\partial D^n_+ \wedge M$. Let $K$ be an object of
$M\dTop$.

\begin{enumerate}
\item We call $K$ {\it cellular\/} if $K$ can be obtained from a point by
  attaching (possibly infinitely many) cells, not necessarily in order of
  increasing dimension.
\item We call $K$ {\it cofibrant\/} if $K$ is a retract of a cellular
  space. The full subcategory of $M\dTop$ consisting of cofibrant spaces is
  denoted $\C (M)$. If $M = S^0$ is the initial pointed monoid, we abbreviate
  this to~$\C$.
\item The space $K$ is called {\it finite\/} if $K$ can be obtained from a
  point by attaching finitely many cells, not necessarily in order of
  increasing dimension. The full subcategory of $M\dTop$ consisting of finite
  spaces is denoted $\Cf(M)$. If $M = S^0$ is the initial pointed monoid, we abbreviate
  this to~$\Cf$.
\item The space $K$ is called {\it homotopy finite\/} if there is a chain (or
  zigzag) of weak equivalences connecting $K$ to an object of~$\Cf(M)$. The
  full subcategory of $M\dTop$ consisting of cofibrant, homotopy finite spaces
  is denoted $\Chf(M)$. If $M = S^0$ is the initial pointed monoid, we abbreviate
  this to~$\Chf$.
\end{enumerate}

In~(4) a {\it weak equivalence\/} is an equivariant map which is a weak
homotopy equivalence on underlying topological spaces.

\begin{remark}
  \label{rem:model_structures}
  The category $M\dTop$ admits a \textsc{Quillen} model structure with
  weak equivalences as above, and fibrations those maps which are
  \textsc{Serre} fibrations on underlying topological spaces. The
  resulting notion of a cofibrant object coincides with the one given above.
\end{remark}

General arguments from model category theory, or a variation on the
\textsc{Whitehead} theorem, imply that a map $X \rTo Y$ of cofibrant
spaces is a weak equivalence if and only if it is a homotopy
equivalence. In particular, a cofibrant space is weakly contractible
if and only if it is contractible.  Since $M$~is discrete, the
forgetful functor restricts to a functor $\C(M) \rTo \C$. In
particular, objects of $M \dTop$ are {\it well-pointed} in the sense
that the inclusion of the base point has the homotopy extension
property.

\medbreak

We will have occasion to use the following standard fact frequently in
the remainder of the paper:

\begin{lemma}
  \label{lem:puppe_sequence}
  Let $f \colon X \rTo Y$ be a map in~$\Top$ such that its homotopy
  cofibre (reduced mapping cone) is contractible. Then the
  reduced suspension $\Sigma f$ of~$f$ is a homotopy equivalence.
\end{lemma}

\begin{proof}
  For any space $V$ we have an exact sequence of pointed sets
  \[ [X,V] \lTo^{f^*}\relax [Y,V] \lTo\relax [\mathrm{hocofibre}\ f,
  V] \lTo\relax [\Sigma X, V] \lTo^{\Sigma f^*}\relax [\Sigma Y, V] \
  ,\] cf.~\cite[Satz~6]{Puppe-Homotopiemengen}, where $[A,B]$ denotes
  the set of (pointed) homotopy classes of maps $A \rTo B$ in~$\Top$.
  Since $\mathrm{hocofibre}\ f \simeq *$ we have $[\mathrm{hocofibre}\
  f, V] = 0$. Hence $f^*$ is monomorphic in the sense of
  \textsc{Puppe} \cite[Footnote~1]{Puppe-Homotopiemengen}, and $\Sigma
  f^*$ is surjective. It follows that $\Sigma f$ has both a left
  homotopy inverse \cite[\S3.1]{Puppe-Homotopiemengen} and a right
  homotopy inverse \cite[\S3.3]{Puppe-Homotopiemengen}.
\end{proof}

\subsection{Barrier cones and projective toric varieties}
\label{subsec:barr-cones-proj}

We will now recall the construction of toric varieties from polytopes.
Standard references are \textsc{Fulton}'s book \cite{Fulton-toric} and
\textsc{Danilov}'s article \cite{Danilov-toric} which contain a wealth of
information on the general theory of toric varieties. More specifically, to
construct varieties from polytopes see \cite[\S1.5, \S3.4]{Fulton-toric} and
\cite[\S5.8, \S11.12]{Danilov-toric}.

\medskip

Let $P \subset \bR^n$ be a lattice polytope (the convex hull of a finite set
of points in~$\bZ^n$) with non-empty interior. Given a non-empty face $F$
of~$P$ we define the {\it barrier cone\/} $C_F$ of~$P$ at~$F$ as the set of
finite linear combination with non-negative real coefficients spanned by the
set $P - F := \{ p - f \,|\, p \in P \hbox{ and } f \in F\}$. Since $C_F$ is a
cone, the integral points in~$C_F$ form a monoid (with respect to the usual
vector sum). By adding a disjoint basepoint, we thus obtain an abelian pointed
monoid $S_F := (C_F \cap \bZ^n)_+$.

For a commutative ring~$R$, let $\tilde R[S_F] = R[S_F]/R[*]$ denote the reduced monoid
ring. If $\emptyset \not= G \subseteq F$ are faces of~$P$, it can be shown
that $S_F$ is obtained from~$S_G$ by inverting a single element. Thus
$\mathrm{Spec}\, \tilde R[S_F]$ is a principal open subset of $\mathrm{Spec}\,
\tilde R[S_G]$. By gluing the affine schemes $U_F := \mathrm{Spec}\, \tilde R[S_F]$
for all non-empty faces of~$P$ we obtain an $R$-scheme~$X_P$, called the
toric variety associated to~$P$.

It can be shown that~$X_P$ is projective. For $R = \bC$ this follows from
\cite[p.~72]{Fulton-toric}. However, the result remains true for arbitrary
commutative rings~$R$. First of all, instead of $P$ we may consider the
dilated polytope $P_D := nP$ without changing the toric variety (note that the
barrier cones of corresponding faces of $P$ and~$P_D$ are the same).  Next,
the polytope $P_D$ defines a Cartier divisor, hence a line bundle, on~$X_P$ as
explained in \cite[page~72]{Fulton-toric} and \cite[\S11.12]{Danilov-toric};
the construction works over any ring, and in fact the resulting line bundle
can explicitly be described as the linearisation of the objects $\mathcal{O}_P
(n)$ to be introduced in Definition~\ref{DefTwistingSheaf} below. Finally,
this line bundle determines a map from $X_P$ to some projective space which
can be shown to be an embedding using Proposition~II.7.2
of~\cite{Hartshorne-AlgGeom}. It remains to see that the hypotheses of that
Proposition are verified, the main point being the surjectivity of ring
homomorphisms from certain polynomial rings to rings of the form $\tilde
R[S_v]$ for $v$ a vertex of~$P_D$. For this it is enough to verify that for
each~$v$ the monoid~$S_v$ is generated by the set of difference vectors $\{p-v
\,|\, p \in P_D \cap \bZ^n \}$. But this is true since $P_D$ is the $n$th
dilation of an $n$-dimensional polytope, cf.~Lemma~VII.3.8 of
\cite{Ewald-CCAG}.

\medbreak

A quasi-coherent sheaf~$\mathcal{F}$ of $\mathcal{O}_{X_P}$-modules gives
rise, by evaluation on the open affine sets $U_F$, to a collection of modules
$\mathcal{F}(U_F)$ over the various rings~$\tilde R[S_F]$, together with
``restriction maps''. Moreover this data completely determines the
sheaf~$\mathcal{F}$.  So we can define a quasi-coherent sheaf as a functor
with values in $R$-modules
$$ M \colon F(P)_0 \rTo R\hbox{-}\mathbf{Mod},\quad F \mapsto M^F $$
(where $F(P)_0$ is the poset of non-empty faces of~$P$)
together with the following data:
\begin{enumerate}
\item For each $F \in F(P)_0$, the module~$M^F$ is equipped with the
  structure of a $\tilde R[S_F]$-module;
\item For each inclusion $G \subseteq F$ in $F(P)_0$, the associated
  map $M^G \rTo M^F$ is $\tilde R[S_G]$-linear;
\item The adjoint
  $ M^G \tensor_{\tilde R[S_G]} \tilde R[S_F] \rTo M^F$
  of the map above is an isomorphism of $\tilde R[S_F]$-modules.
\end{enumerate}

\subsection{Non-linear sheaves}
\label{subsec:non-linear-sheaves}

\begin{definition}
\label{Def:non-lin-sheaf}
  A {\it non-linear sheaf\/} on~$X_P$ is a functor
  $$ Y \colon F(P)_0 \rTo \Top,\quad F \mapsto Y^F $$
  together with the following data:
  \begin{enumerate}
  \item For each $F \in F(P)_0$, the space~$Y^F$ is equipped with a base point
    preserving (right) action of the pointed monoid~$S_F\,$;
  \item For each inclusion $G \subseteq F$ in $F(P)_0$, the associated map
    $Y^G \rTo Y^F$ is $S_G$-equivariant;
  \item Using the notation of~(2), let $Y^G_c \rTo^g Y^G$ be a cofibrant
    replacement of~$Y^G$, \ie, $Y^G_c \in \C(S_G)$ and $g$ is an
    $S_G$-equivariant weak homotopy equivalence. Then the map
    $$ Y^G_c \wedge_{S_G} S_F \rTo Y^F\ ,$$
    adjoint to the composition $Y^G_c
    \rTo^\sim Y^G \rTo Y^F$, is a weak equivalence.
  \end{enumerate}
\end{definition}

Existence of cofibrant replacements is a direct consequence of the
model category structure mentioned in
Remark~\ref{rem:model_structures}. Using a cofibrant replacement
ensures that the ``gluing condition''~(3) is weakly homotopy
invariant. Moreover, standard model category arguments show that we
could equivalently have worked with a fixed cofibrant replacement, or
we could have asked for the gluing condition to be satisfied for all
cofibrant replacements.  In particular, we obtain:

\begin{lemma}
\label{lem:sheaf-cond}
  Suppose $Y \colon F(P)_0 \rTo \Top$ is a diagram satisfying conditions~(1)
  and~(2) above. Suppose moreover that $Y$ is locally cofibrant in the
  sense that $Y^F \in \C(S_F)$ for all $F \in F(P)_0$. Then $Y$ satisfies the
  gluing condition~(3) if and only if for all inclusions $G \subseteq F$
  in $F(P)_0$, the map
  $$ Y^G \wedge_{S_G} S_F \rTo Y^F \ ,$$
  adjoint to the structure map $Y^G \rTo Y^F$, is a weak equivalence.
  \doqed
\end{lemma}

\begin{definition}
\label{def:sheaves-cofibrancy}
  \begin{enumerate}
  \item A non-linear sheaf $Y$ on $X_P$ is called {\it weakly cofibrant\/} if
    for all $F \in F(P)_0$ the component~$Y^F$ is cofibrant as a pointed
    topological space, \ie, $Y^F \in \C$. The category of weakly cofibrant non-linear sheaves
    on~$X_P$ is denoted~$\wchcoh$.
  \item A non-linear sheaf $Y$ on $X_P$ is called {\it locally cofibrant\/} if
    $Y^F \in \C(S_F)$ for all $F \in F(P)_0$. The category of locally
    cofibrant non-linear sheaves on~$X_P$ is denoted~$\hcoh$.
  \item A map of non-linear sheaves is called a {\it weak equivalence\/} if
    all its components are weak homotopy equivalences of spaces.
  \end{enumerate}
\end{definition}

The notation $\hcoh$ is intended to suggest that a non-linear sheaf is a
homotopy-theoretic version of a quasi-coherent sheaf. Every locally cofibrant
non-linear sheaf is weakly cofibrant.

\medbreak

The most important examples of non-linear sheaves are the ``twisting sheaves'',
formed by using translates of the monoids~$S_F$.

\begin{figure}[ht]
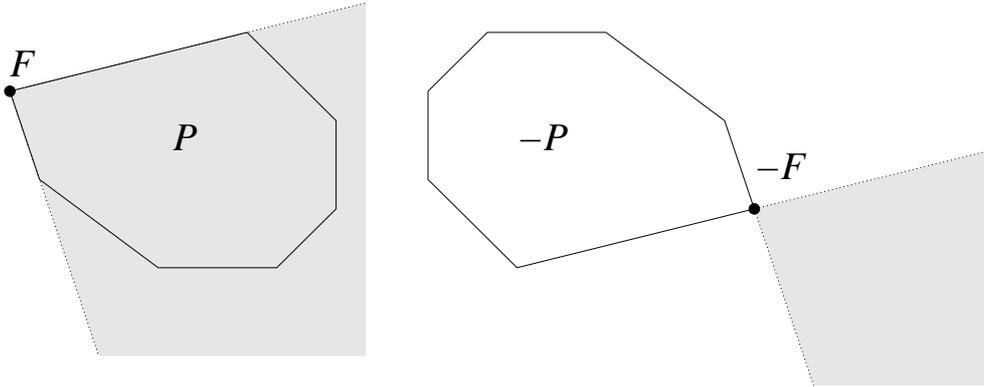

\begin{center}
\input barrier.pstex_t
\caption{The construction of~$\mathcal{O}_P(k)$}
\label{fig:barriercones}
\end{center}
\end{figure}

\begin{definition}
\label{DefTwistingSheaf}
  For all $k \in \bZ$, we define the $k$th {\it twisting sheaf\/}, denoted
  $\mathcal{O}_P (k)$, as the non-linear sheaf
  $$ \mathcal{O}_P (k) \colon F(P)_0 \rTo \Top,\quad
     F \mapsto \Big( \bZ^n \cap (C_F + kF) \Big)_+ $$
  where $C_F$ is the barrier cone of~$P$ at~$F$, and
  $$ C_F + kF = \{ x+kf \,|\, x \in C_F \hbox{ and } f \in F \} \ . $$
\end{definition}
Note that $\mathcal{O}_P (0)^F = S_F$ and $\mathcal{O}_P (k)^F \iso S_F$ (not
canonically). Moreover, $\mathcal{O}_P(k)$ is a non-linear sheaf by
Lemma~\ref{lem:sheaf-cond}. See Figure~\ref{fig:barriercones} for a picture of
$\mathcal{O}_P(1)^F$ and $\mathcal{O}_P(-1)^F$ (the shaded areas) for $F$ a
vertex of~$P$.

By passage to reduced free modules, we obtain a diagram $F \mapsto \tilde R
[\mathcal{O}_P (k)^F] $ which is a quasi-coherent sheaf in the sense of
\S\ref{subsec:barr-cones-proj}; this is the algebraic geometers' $k$th
twisting sheaf on~$X_P$.

\begin{definition}
  \label{def:twisting}
  For $Y,Z \in \hcoh$ we define their {\it tensor product\/} $Y \tensor Z$ by
  $$ Y \tensor Z \colon F(P)_0 \rTo \Top,\quad F \mapsto Y^F \wedge_{S_F} Z^F
  \ .$$
  Here $Y^F \wedge_{S_F} Z^F$ is the co-equaliser of the two maps $Y^F \wedge
  S_F \wedge Z^F \rTo Y^F \wedge Z^F$ given by the action of~$S_F$ on~$Y^F$
  and~$Z^F$, respectively.

  For $j \in \bZ$ we define the $j$th {\it twist\/} of~$Y$ as $Y(j) := Y
  \tensor \mathcal{O}_P(j)$. Both $Y \tensor Z$ and $Y(j)$ are objects
  of~$\hcoh$ again. The twisting functor $Y \mapsto Y(j)$ will also be
  denoted $\theta_j \colon \hcoh \rTo \hcoh $.
\end{definition}

Note that the isomorphism $\mathcal{O}_P (k)^F \iso S_F$ induces a non-canonical
isomorphism $Y(j)^F \iso Y^F$.  It is easy to check that $\mathcal{O}_P (j)
\tensor \mathcal{O}_P (k) \iso \mathcal{O}_P (j+k)$ and thus $Y(j)(k) \iso
Y(j+k)$.  Moreover, $Y(0) \iso Y$, so twisting defines a self-equivalence of
the category~$\hcoh$ which maps weak equivalences to weak equivalences.

\begin{definition}
\label{def:canonical_sheaves}
  Given a space $K \in \Top$ we define the diagram $K \wedge \mathcal{O}_P
  (k)$ by $(K \wedge \mathcal{O}_P(k))^F := K \wedge \mathcal{O}_P(k)^F$.  The
  functor
  $$ \psi_k \colon \C \rTo \hcoh,\quad K \mapsto K \wedge \mathcal{O}_P (k) $$
  is called the $k$th {\it canonical sheaf\/} functor.
\end{definition}

From the remarks above we have isomorphisms $\theta_j \circ \psi_k (K) \iso \psi_{k+j}
(K)$ which are natural in~$K$.

\subsection{Total cofibres}
\label{subsec:total-cofibres}

Let $P \subset \bR^n$ be a polytope.  The set $F(P)$ of all faces of $P$ is
partially ordered by inclusion and can thus be considered as a category with
initial object $\emptyset$ and terminal object $P$. We define $F(P)^1 = F(P)
\setminus \{P\}$.

\begin{definition}
\label{DefTotalCofibre}
  Given a functor
  $\, Y \colon F(P) \rTo \C,\ F \mapsto Y^F\, $
  we define the {\it total cofibre\/} of~$Y$, denoted $\Gamma(Y)$, as the
  cofibre of the canonical cofibration
  $$ \mathrm{hocolim}\, Y|_{F(P)^1} \rTo \mathrm{hocolim}\, Y \ . \eqno {(*)} $$
\end{definition}

\begin{remark}
  \label{rem:def-total-cofibres}
  \begin{enumerate}
  \item Since $F(P)$ has the terminal object~$P$ the space
    $\mathrm{hocolim}\, Y$ is weakly homotopy equivalent to $Y^P$, and the
    total cofibre of~$P$ is weakly equivalent to the homotopy cofibre of the
    map $\mathrm{hocolim}\, Y|_{F(P)^1} \rTo Y^P$.
  \item The definition of $\Gamma (Y)$ depends on the combinatorial type
    of~$P$ only, not on its actual geometry. If $P = \Delta^{n-1}$ is a
    simplex, this definition coincides with the usual definition of the total
    cofibre of an $n$-cubical diagram as given in
    \cite[Definition~1.4]{Goodwillie-Calculus2}.
  \end{enumerate}
\end{remark}

A diagram $Y \colon F(P)_0 \rTo \C$, eg., a weakly cofibrant
non-linear sheaf, can be considered as a diagram defined on all
of~$F(P)$ by setting $Y^\emptyset = *$, so $\Gamma(Y)$ is defined in
this case as well.

\subsubsection*{Iterating the total cofibre construction}

\begin{lemma}
\label{lem:GammaGamma}
Suppose $P$ and~$Q$ are polytopes, and suppose~$Z$ is a diagram
$ Z \colon F(P) \times F(Q) \rTo \C,\ (F,G) \mapsto Z_F^G$.
  There is a natural homeomorphism
  $$ \Gamma \left( F \mapsto \Gamma (Z_F^?) \right) \ \iso
     \ \Gamma \left( G \mapsto \Gamma (Z_?^G) \right) \ . $$
\end{lemma}

\medbreak
\vglue 1 \smallskipamount
\begin{proof}
  The proof is encoded into the following diagram:
  \begin{diagram}
    \mathop{\mathrm{hocolim}}_{F \in F(P)^1}\,
      \mathop{\mathrm{hocolim}}_{G \in F(Q)^1}\, Z_F^G & \rTo[l>=3em] &
      \mathop{\mathrm{hocolim}}_{F \in F(P)}\,
      \mathop{\mathrm{hocolim}}_{G \in F(Q)^1}\, Z_F^G & \rTo[l>=3em] &
      \mathop{\mathrm{hocolim}}_{G \in F(Q)^1}\, \Gamma (Z_?^G) \\
    \dTo && \dTo && \dTo \\
    \mathop{\mathrm{hocolim}}_{F \in F(P)^1}\,
      \mathop{\mathrm{hocolim}}_{G \in F(Q)}\, Z_F^G & \rTo &
      \mathop{\mathrm{hocolim}}_{F \in F(P)}\,
      \mathop{\mathrm{hocolim}}_{G \in F(Q)}\, Z_F^G & \rTo &
      \mathop{\mathrm{hocolim}}_{G \in F(Q)}\, \Gamma (Z_?^G) \\
    \dTo && \dTo && \dTo \\
    \mathop{\mathrm{hocolim}}_{F \in F(P)^1}\, \Gamma (Z_F^?) & \rTo &
      \mathop{\mathrm{hocolim}}_{F \in F(P)}\, \Gamma (Z_F^?) & \rTo &
      \Gamma \Gamma (Z_?^?) \\
  \end{diagram}
  All rows and columns are cofibre sequences. For example, the first
  column is obtained by applying the functor $\mathrm{hocolim}_{F \in F(P)^1}$
  to the cofibre sequence defining $\Gamma (Z_F^?)$, and the first row is
  obtained by applying the functor $\mathrm{hocolim}_{G \in F(Q)^1}$ to the
  cofibre sequence defining $\Gamma (Z_?^G)$; note that homotopy colimits
  commute among themselves as well as with taking cofibres.
\end{proof}

\subsubsection*{A vanishing theorem for total cofibres}

\begin{observation}
  An $(n+1)$-cubical diagram $Y \colon F(\Delta^n) \rTo \C$ can be written
  as a map $f \colon Z_1 \rTo Z_2$ of $n$-cubical diagrams: If $v$ is a vertex
  of~$\Delta^n$, then $Z_1$ is the restriction of~$Y$ to the poset of all
  faces of~$\Delta^n$ not containing~$v$, and $Z_2$ is the restriction of~$Y$
  to the poset of all faces of~$\Delta^n$ containing~$v$. The components
  of~$f$ are the structure maps $Y^F \rTo Y^{F \vee \{v\}}$ of~$Y$ for $v
  \notin F \in F(\Delta^n)$. (For $n=1$, the diagram $Y$ is a square, and $f$
  is the map from the top to the bottom arrow, or the map from the left to the
  right arrow.)  If $f$~consists of weak equivalences, the diagram~$Y$ is
  homotopy cocartesian (see the remarks preceding Definition~1.4
  in~\cite{Goodwillie-Calculus2}), and its total cofibre is weakly
  contractible.
\end{observation}

The point is that $\Gamma (Y)$ is homeomorphic to the total cofibre of the
$n$-cubical diagram $\mathrm{hocofibre}\,f$. If $f$ consists of weak
equivalences, $\mathrm{hocofibre}\,f$ consists of weakly contractible spaces
only, so its total cofibre is homotopically trivial.

\medskip

We will prove the following generalisation of this simple vanishing
criterion (an essential ingredient for the proofs of
Lemma~\ref{lem:sigma-is-suspension} and
Proposition~\ref{prop:psi0Gamma-is-suspension}):

\begin{theorem}
\label{GammaVanishes}
Suppose the functor $ Y \colon F(P) \rTo \C $ has the property that
for some non-empty face~$A$ of~$P$ ``all structure maps in
$A$-direction are weak equivalences'', \ie, for all $ F \in F(P) $ the
map $ Y^F \rTo Y^{F\vee A}$ is a weak equivalence.  Then the total
cofibre of~$Y$ is weakly contractible.
\end{theorem}

\bigbreak

The total cofibre~$\Gamma(Y)$ measures the deviation of~$Y$ from being
a homotopy colimit diagram. If $Y^P \simeq \hocolim (Y|_{F(P)^1})$,
\ie, if the canonical map~$(*)$ of Definition~\ref{DefTotalCofibre} is
a (weak) homotopy equivalence, then $\Gamma(Y)$ is (weakly)
contractible. Conversely, if $\Gamma (Y)$ is contractible, the
canonical map~$(*)$ suspends to a weak equivalence
(Lemma~\ref{lem:puppe_sequence} and
Remark~\ref{rem:def-total-cofibres}~(1)). In this sense, a vanishing
result for total cofibres is nothing but a ``computation'' of a
homotopy colimit, up to suspension.

\bigbreak

We begin with some technical preliminaries.
To simplify the notation, we define three subcategories of~$F(P)^1$:
\begin{eqnarray}
        \mathcal{C}_0 &:=& \{ \emptyset \} \cup lk (A) \nonumber \\
        \mathcal{C}_1 &:=& \{ \emptyset \} \cup  \overline {st} (A) \nonumber \\
        \mathcal{C}_2 &:=& \{ \emptyset \} \cup \overline {ast} (A) \nonumber
\end{eqnarray}
Links, stars and antistars are computed in the complex~$F(P)_0^1$ unless
indicated otherwise. Note that $\mathcal{C}_1 \cap \mathcal{C}_2 = \mathcal{C}_0$.

Let $ \iota \colon st(A) \rTo \mathcal{C}_1 $ denote the inclusion,
and define
\[\Phi \colon \mathcal{C}_1 \rTo st(A), \quad F \mapsto F \vee A\]
(this is well defined since $ F \vee A \not= P $ by~\ref{CombDesc}~(2)).
Then $\Phi \circ \iota =\id_{st(A)} $, and there is a natural
transformation of functors $\theta \colon \id \rTo \iota \circ \Phi$
with $F$-component given by the inclusion $ F \rTo F \vee A $.

\begin{lemma}
  \label{iota}
  The inclusion $ \iota \colon st(A) \rTo \mathcal{C}_1 $ induces a
  homotopy equivalence $ \iota_* \colon \hocolim\ Y|_{st(A)} \rTo
  \hocolim\ Y|_{\mathcal{C}_1} $ with homotopy inverse~$\gamma$
  induced by~$\Phi$ and~$\theta$.
\end{lemma}

\begin{proof}
  This follows from \cite[Corollary~3.14]{WZZ-hocolims}. We provide a
  proof for the reader's convenience.  The map $\gamma$ induced
  by~$\Phi$ and~$\theta$ factors as
  \[\hocolim\ Y|_{\mathcal{C}_1} \rTo^{\theta_*} \hocolim\ \big(
  \Phi^* (Y|_{st(A)})
  \big) \rTo^{\Phi_*} \hocolim\ Y|_{st(A)}\]
  where the first map is induced by the natural
  transformation~$\theta$, and the second map is induced
  by~$\Phi$. The composition $\Phi_* \circ \theta_* \circ \iota_* =
  \gamma \circ \iota_*$
  is the identity map of $\hocolim\
  Y|_{st(A)}$ since $\Phi \circ \iota = \id$.

  We are left to show that $\iota_* \circ \gamma$ is homotopic to the
  identity map of $\mathrm{hocolim}\ Y|_{\mathcal{C}_1}$. The natural
  transformation $\theta \colon \id \rTo \iota \circ \Phi$ can be
  encoded as a single functor
  \[\Upsilon \colon \mathcal{C}_1 \times [1] \rTo \mathcal{C}_1\]
  (where $[1] = \{0 \rTo 1\}$ is the category with two objects and a
  single non-trivial morphism) such that $\Upsilon|_{\mathcal{C}_1
    \times \{0\}} = \id$ and $\Upsilon|_{\mathcal{C}_1 \times \{1\}} =
    \iota \circ \Phi$. Now $\mathrm{hocolim}\ \Upsilon^* (Y|_{\mathcal{C}_1})$ is
  homeomorphic to the mapping cylinder~$Z_\gamma$ of the map~$\gamma$,
  and the functor $\Upsilon$ induces a map
  \[Z_\gamma \iso \mathrm{hocolim}\ \Upsilon^* (Y|_{\mathcal{C}_1}) \rTo
  \mathrm{hocolim}\ Y|_{\mathcal{C}_1} \ .\]
  Pre-composition with the map $\left(\mathrm{hocolim}\
    Y|_{\mathcal{C}_1}\right) \times [0,1] \rTo Z_\gamma$ yields the
  desired homotopy.
\end{proof}

\begin{lemma}
\label{Psi}
  Let~$\Psi$ denote the composition $ \mathcal{C}_0 \rTo \mathcal{C}_1 \rTo^\Phi st(A) $.
  Then~$\Psi$ induces a homotopy equivalence $ \alpha \colon \hocolim\ \Psi^* (Y|_{st(A)})
  \rTo \hocolim\ Y|_{st(A)} $.
\end{lemma}

\begin{proof}%
It suffices to show that~$\Psi$ is right cofinal
\cite[dual of Theorem~XI.9.2]{BK-Monster}\cite[Proposition~3.10]{WZZ-hocolims}, 
\ie, for all elements $ B \in st(A) $ the category $ B \downarrow \Psi $ is
contractible.

Case~1: $B = A$. Then $ B \downarrow \Psi = \mathcal{C}_0 $ has the initial
object~$\emptyset$, hence has contractible classifying space.

Case~2: $B \supset A$.
Then $ B \downarrow \Psi = \{ F \in lk(A) | B \subseteq F \vee A \} $
by definition of~$\Psi$. We also have the equality
$$ B \downarrow \Psi = \{ F \in lk(B) | B \subseteq F \vee A \} \ .$$
Indeed, using~\ref{CombDesc}~(2), we conclude that for every $F \in B \downarrow \Psi$
we have $ B \not\subseteq F$ since~$F$ does not contain~$A$, and
$$ F \vee B = F \vee A \vee B = F \vee A \not= P $$
(since $ B \subseteq F \vee A $ by definition of $ B \downarrow \Psi $),
thus $ F \in lk(B) $ by Lemma~\ref{CombDesc}~(2).
Conversely, if $F \in lk (B)$ satisfies $ B
\subseteq F \vee A$, we have $ A \not\subseteq F $ since otherwise
$ B \subseteq F \vee A = F $ which contradicts $ F \in lk(B) $. Moreover,
$ F \vee A \subseteq F \vee B \not= P $, and we conclude
$ F \in B \downarrow \Psi $.

By Corollary~\ref{CombStarLinkCor} we know
\[lk (B) \setminus (B \downarrow \Psi) = 
\{ F \in lk(B) | B \not\subseteq F \vee A \} =
\overline {st}_{lk (B)} (A)\]
and consequently $ B \downarrow \Psi = ast_{lk(B)}(A)$.
Now $|lk(B)| \pliso S^{n-2} $ since the boundary of~$P$ is \hbox{$PL$-homeomorphic}
to an $(n-1)$-sphere. Thus $|\overline{st}_{lk(B)}(A)|$ is an
$(n-2)$-dimensional ball. We can now apply Corollary~\ref{RegularNbhdSpheres} to show
that the classifying space of $B \downarrow \Psi$ is contractible.
\end{proof}

\begin{lemma}
\label{Xi}
The inclusion $ \Xi \colon ast(A) \rTo \mathcal{C}_2 $ induces a homotopy equivalence
$ \delta \colon \hocolim\ \Xi^* (Y|_{\mathcal{C}_2}) \rTo \hocolim\ Y|_{\mathcal{C}_2} $.
\end{lemma}

\begin{proof}%
  It suffices to show that~$\Xi$ is right cofinal
  \cite[Theorem~XI.9.2]{BK-Monster}\cite[Proposition~3.10]{WZZ-hocolims},
  \ie, for all elements $ G \in \mathcal{C}_2 $ the category $ G
  \downarrow \Xi $ is contractible.  Fix an object $ G \in
  \mathcal{C}_2 $.

Case~1: $ G \in ast(A)$. Then $ G \downarrow \Xi $ contains~$G$ as an
initial object and hence is contractible.

Case~2: $G = \emptyset$. Then $ G \downarrow \Xi = ast(A) $. Its classifying space
is contractible by Corollary~\ref{RegularNbhdSpheres}, applied to $K =
F(P)_0^1$ and $ C = ast (A) $.

Case~3: $ G \in lk(A) $. Then $ G \downarrow \Xi = \{ F \in ast(A) | G \leq F
\} = st(G) \cap ast(A) $; this is an order filter in $F(P)_0^1$. Its
complement is $Z:= \overline{ast}(G) \cup \overline{st}(A)$.
Now
$$ \overline{ast}(G) \cap \overline {st}(A) = \overline {st}(A) \setminus st(G)
   = \overline {st}(A) \setminus st_{\overline {st}(A)}(G) \ .$$
Consequently, we can write
$$ Z = \overline{ast}(G) \cup \overline{st}(A)
   = \overline{ast}(G) \amalg st_{\overline{st}(A)}(G)
   = \overline{ast}(G) \cup \overline{st}_{\overline{st}(A)}(G) $$
(where the last equality holds since $Z$~is a complex anyway, thus using the
closed star instead of the open star does not make any difference).
Thus $|Z|$~is the union of the two $(n-1)$-dimensional balls $|\overline{ast}(G)|$
and~$|\overline{st}_{\overline{st}(A)}(G)|$; their intersection is
$lk_{\overline{st}(A)}(G)$ which is an $(n-2)$-dimensional
ball since $G \in lk(A)$ (whence \hbox{$\hat G \in \partial |\overline{st}(A)|\,$}).
We conclude that $|Z|$~is an $(n-1)$-dimensional ball \cite[Corollary~3.16]{RS-PL}.
Now Corollary~\ref{RegularNbhdSpheres}, applied to $K = F(P)_0^1$ and $C=st(G)
\cap ast(A)$, shows that $|N(G\downarrow \Xi)| = |N (st(G)\cap ast(A))| \simeq
*$ as claimed.
\end{proof}

\bigbreak
\noindent {\bf Proof of Theorem~\ref{GammaVanishes}.}
Since the categories~$\mathcal{C}_1$ and~$\mathcal{C}_2$ form a convex cover
\cite[\S0]{Goodwillie-Calculus2} of~$F(P)^1$ with intersection~$\mathcal{C}_0$,
Proposition 0.2, {\it op.\thinspace{}cit.}, shows that the following square is homotopy cocartesian:
\begin{diagram}[l>=3em]
        \hocolim\ Y|_{\mathcal{C}_0} & \rTo & \hocolim\ Y|_{\mathcal{C}_1} \\
        \dTo && \dTo \\
        \hocolim\ Y|_{\mathcal{C}_2} & \rTo & \hocolim\ Y|_{F(P)^1} \\
\end{diagram}
In particular, the space~$\Gamma (Y)$ is weakly homotopy equivalent to the
total cofibre of the following square (we have used
Remark~\ref{rem:def-total-cofibres}~(1) to replace $\mathrm{hocolim}\, Y$ by
the weakly equivalent space $Y^P$ in the terminal entry):
\begin{diagram}[l>=3em,eqno=(*)]
        \hocolim\ Y|_{\mathcal{C}_0} & \rTo^f & \hocolim\ Y|_{\mathcal{C}_1} \\
        \dTo && \dTo \\
        \hocolim\ Y|_{\mathcal{C}_2} & \rTo_g & Y^P \\
\end{diagram}
We will show that~$f$ and~$g$ are weak homotopy equivalences. Then their
homotopy cofibres are weakly contractible, and since
$$ \Gamma (Y) \iso \mathrm {hocofibre}\ (\mathrm {hocofibre}\ (f) \rTo \mathrm
{hocofibre}\ (g)) $$ this proves the assertion of the theorem.

\medbreak

We can embed the square~$(*)$ into
the bigger commutative diagram shown in Fig.~\ref{fig:diagram}.
\begin{figure}[htbp]
\begin{diagram}[l>=3em,PS]
  \hocolim\ \Psi^* (Y|_{st(A)}) & \rTo^\alpha & \hocolim\ Y|_{st (A)} \\
  \uTo<\beta && \uTo>\gamma \\
  \hocolim\ Y|_{\mathcal{C}_0} & \rTo^f & \hocolim\ Y|_{\mathcal{C}_1} \\
  \dTo && \dTo \\
  \hocolim\ Y|_{\mathcal{C}_2} & \rTo^g & Y^P \\
  \uTo<\delta && \uTo>\epsilon & \luTo>{\mathrm {pr}_2} \\
  \hocolim\ \Xi^*(Y|_{\mathcal {C}_2}) & \rTo_\eta & \hocolim\ Y(P)^{ast(A)}
  & \rTo^{\iso}& |Nast(A)| \times Y^P \\
\end{diagram}
\caption{The diagram used in the proof of Theorem~\ref{GammaVanishes}}
\label{fig:diagram}
\end{figure}
(Here $Y(P)^{ast(A)}$ denotes the constant diagram on~$ast(A)$ with
value~$Y^P$.)  The map~$\alpha$ is induced by~$\Psi$; it is a weak equivalence
by Lemma~\ref{Psi}. Similarly, $\Phi$~induces the weak equivalence~$\gamma$ by
Lemma~\ref{iota}. The map~$\beta$ is induced by the natural transformation $
Y|_{\mathcal {C}_0} \rTo \Psi^* (Y|_{st(A)}) $ with $F$-components given by $
Y^F \rTo Y^{F \vee A} $. But the latter are weak homotopy equivalences by
hypothesis on~$Y$. Hence~$\beta$ is a weak equivalence by the Homotopy Lemma
\cite[Lemma~XII.4.2]{BK-Monster} \cite[Lemma~4.6]{WZZ-hocolims}. This proves
that~$f$ is a weak equivalence as well.

Since~$|F(P)_0^1| \pliso S^{n-1}$ application of
Corollary~\ref{RegularNbhdSpheres} yields \hbox{$|Nast(A)| \simeq *$}.  It
follows that $\epsilon$~is a weak equivalence. The map~$\delta$, induced
by~$\Xi$, is a weak equivalence by Lemma~\ref{Xi}. Finally, the map~$\eta$ is
induced by the natural transformation $\Xi^* (Y|_{\mathcal{C}_2}) \rTo
Y(P)^{ast(A)}$ with $F$-components given by the weak homotopy equivalences
\hbox {$Y^F \rTo Y^{F \vee A}$} (recall from~\ref{CombDesc}~(3) that $F \vee A
= P$ for all $F \in ast(A)$). Hence $\eta$~is a weak equivalence itself by the
Homotopy Lemma \cite[Lemma~XII.4.2]{BK-Monster}
\cite[Lemma~4.6]{WZZ-hocolims}.  This proves that $g$~is a weak homotopy
equivalence as claimed.  \doqed

\subsection{Total cofibres of canonical sheaves}
\label{subsec:total-cofib-can}

Let $P \subset \bR^n$ be a lattice polytope with non-empty interior. For any
integer $k \in \bZ$ we define $kP := \{ kp \,|\, p \in P\}$.

\begin{theorem}
  \label{thm:Ehrhart}
  {\rm \cite[\S{}IV.6]{Ewald-CCAG}}
  There is a polynomial $E_P (T) \in \bQ[T]$ of degree~$n$ with the following
  properties:
  \begin{enumerate}
  \item If $k \geq 0$ is an integer, then $E_P(k) = \# \big( kP \cap \bZ^n
    \big)$. In particular, $E_P(0)=1$.
  \item If $k < 0$ is an integer, then $(-1)^n E_P (k) = \# \big(
    \mathrm{int}\, (kP) \cap \bZ^n \big)$.\doqed
  \end{enumerate}
\end{theorem}

The polynomial $E_P(T)$ of the theorem is called the {\it \textsc{Ehrhart\/}
polynomial\/} of~$P$.

\medbreak

For a non-empty proper face $F$ of~$P$ let $T_F$ denote the supporting cone
of~$F$; it is the intersection of all supporting half-spaces containing $F$ in
their boundary. (Of course it is enough to restrict to facet-defining
half-spaces.) By convention $T_P = \bR^n$.

Let~$C_F$ denote the barrier cone (\S\ref{subsec:barr-cones-proj}) of~$P$
at~$F$; it is the set of linear combinations with non-negative real
coefficients spanned by $P-F = \{ p-f | p \in P,\ f \in F\} $.  Using \textsc
{Farkas\/}' lemma (\cite[\S1.4]{Ziegler-Polytopes} or
\cite[Lemma~I.3.5]{Ewald-CCAG}) it can be shown that $F + C_F = T_F$.
Moreover, every polytope is the intersection of all its supporting
half-spaces, thus $P = \bigcap_{F \in F(P)_0^1} T_F$.

\bigbreak

Recall that for~$k \in \bZ$ the twisting sheaf $\mathcal{O}_P(k)$ is defined as
$$ \mathcal{O}_P(k) \colon F(P)_0 \rTo \Top, \quad
F \mapsto \left((kF + C_F) \cap \bZ^n\right)_+ \ ,$$
and for $K \in \Top$ the $k$th canonical sheaf is defined as $\psi_k (K) = K
\wedge \mathcal{O}_P (k)$.  Note that $\mathcal{O}_P(1)^F = (T_F \cap
\bZ^n)_+$ and $\mathcal{O}_P(0)^F = (C_F \cap \bZ^n)_+ = S_F$.

\bigbreak

The following theorem generalises
\cite{H-Proj}, Corollaries~3.7.4--5 (the case $P = \Delta^n =
\mathrm {conv}\ \{0,\,e_1,\,\ldots,\,e_n\}$).

\begin{theorem}
\label{HocolimT}
  Suppose $P \subseteq \bR^n$ is a lattice polytope with
  non-empty interior. Let $K \in \C$ be a cofibrant pointed topological space.
  \begin{enumerate}
  \item For every integer $k \geq 0$ there is a natural homotopy equivalence
    $$ \Gamma \left( \psi_k (K) \right) \simeq \left( kP \cap \bZ^n
    \right)_+ \wedge S^n \wedge K = \bigvee_{E_P (k)} \Sigma^n K \ . $$
  \item For every integer $k < 0$ there is a natural homotopy equivalence
    $$ \Gamma \left( \psi_k (K) \right) \simeq \left( (\mathrm{int}\, kP)
    \cap \bZ^n \right)_+ \wedge K = \bigvee_{(-1)^n E_P (k)} K \ . $$
    In particular, $\Gamma \left( \psi_k(K) \right) \simeq *$ if the interior of
    $kP$ does not contain any lattice point (\ie, if $E_P (k) = 0$).
  \end{enumerate}
\end{theorem}

\begin{proof}%
Since homotopy colimits commute with smash products there is a canonical
isomorphism $\Gamma (\psi_k(K)) = \Gamma (\psi_k(S^0 \wedge K)) \iso \Gamma
(\psi_k(S^0)) \wedge K$. It is thus enough to consider the case $K = S^0$
only. Note that $\psi_k (S^0) \iso \mathcal{O}_P (k)$.
The space $ \Gamma \left( \mathcal{O}_P(k) \right) $ is homeomorphic to the homotopy
cofibre of the natural map
$$ \kappa \colon \mathrm {hocolim}\, \mathcal{O}_P(k)|_{F(P)_0^1} \rTo (\bZ^n)_+ $$
which is induced by the inclusions $ k F + C_F \subseteq k P + C_P = \bR^n $.
Define, for fixed~$x \in \bZ^n$, the functor with values in (unpointed) topological spaces
$$ T(k)^x \colon F(P)_0^1 \rTo \mathbf{Top}, \quad F \mapsto \{x\} \cap
\mathcal{O}_P(k)^F \ .$$
There is a natural isomorphism of functors $ \mathcal{O}_P(k)|_{F(P)_0^1}
\iso \left(\amalg_{x \in \bZ^n} T(k)^x\right)_+$.  Consequently, there is a
homeomorphism
$$ \mathrm {hocolim}\, (\mathcal{O}_P(k)|_{F(P)_0^1}) \iso \left(
\amalg_{x \in \bZ^n} \mathrm {hocolim}^\prime\, T(k)^x \right)_+ $$
where $\mathrm{hocolim}^\prime$ denotes the {\it unpointed\/} homotopy
colimit. This homeomorphism induces, for any point $x \in \bZ^n$, a
homeomorphism $ \kappa^{-1} (x) \iso \mathrm {hocolim}\, T(k)^x $.

\bigbreak

To prove~(1) it is thus sufficient to show that $\mathrm {hocolim}^\prime\,
T(k)^x \iso S^{n-1}$ for \hbox{$x \in \bZ^n \cap kP$} and $\mathrm
{hocolim}^\prime\, T(k)^x \simeq *$ for $x \notin \bZ^n \cap kP $.

\medbreak

Assume $ k > 0 $ first. It is enough to consider the case $k=1$ since the
functors $\mathcal{O}_P (k)$ and $\mathcal{O}_{kP}(1)$ are isomorphic, and we
have an equality $E_P (kT) = E_{kP} (T)$. So assume $k=1$. Then
$\mathcal{O}_P(1)^F = (T_F \cap \bZ^n)_+$ is the intersection of the
supporting cone of~$F$ with~$\bZ^n$ (plus a disjoint base point). Fix a point
$x \in \bZ^n$.

If $x \in P$, the functor $T(k)^x$~is the constant functor with a
one point space as value, hence $\mathrm {hocolim}^\prime\, T(k)^x \iso |NF(P)_0^1|
\iso S^{n-1}$ by Lemma~\ref{Triangulation}.

Now assume $ x \notin P$. Let~$F$ denote a proper non-empty face of~$P$. From
Lemma~\ref{VisibleFacets} and the definition of supporting cones we conclude that
$x \in T_F$ if and only if $F$~is invisible from~$x$.
In particular, $T(k)^x (F) = \{x\}$ if $F \in \Invis(x)$, and $T(k)^x (F) =
\emptyset$ if $F \notin \Invis(x)$. By definition of homotopy colimits, this
implies $$ \mathrm{hocolim}^\prime\, T(k)^x \iso |N\Invis(x)| \ .$$
But by Corollary~\ref{BIcontractible}, this space is contractible.

\medbreak

Now assume $k = 0$. Then $\mathcal{O}_P(0)^F = (C_F \cap \bZ^n)_+$ is the
intersection of the barrier cone of~$P$ at~$F$ with~$\bZ^n$ (plus a disjoint
base point). Fix a point $x \in \bZ^n$.

If $x=0$ we see $$\hocolim^\prime\ T(0)^0 \iso |NF(P)_0^1| \pliso S^{n-1} $$
since by their definition all barrier cones contain the origin, \ie,
$T(0)^0$~is the constant functor with value~$\{0\}$ in this case.

If $x \not= 0$, let $N_F := C_F^{\vee}$ denote the dual cone of~$C_F$; it is given by
$$ N_F = \left\{ v \in \bR^n \,|\, \forall p \in C_F \colon
   \langle p, \, v \rangle \geq 0 \right\} \ .$$
It can be shown that $N_F$~is the cone of inward pointing normal vectors
of~$F$, and that the dual of~$N_F$, given by
$$ N_F^\vee := \left\{ p \in \bR^n \,|\, \forall v \in N_F \colon
   \langle v, \, p \rangle \geq 0 \right\} \ , $$
is the barrier cone~$C_F$ \cite[\S{}I.4 and~\S{}V.2]{Ewald-CCAG}.

Let $U(x)$ denote the poset of all non-empty proper faces~$F$ of~$P$ satisfying $x \in
\mathcal{O}_P(0)^F = (C_F \cap \bZ^n)_+ $. Then $ \hocolim^\prime\ T(0)^x \iso |NU(x)| $.
By the above we have equivalences
$$ F \in U(x) \iff x \in C_F = N_F^\vee
   \iff \forall v \in N_F \colon \langle -x, \, v \rangle \leq 0 \ . $$
This means that $U(x) = \Up(-x)$ is the set of upper faces of~$P$
(with respect to~$-x$) in the sense
of~Definition~\ref{LowUpDef}. By Corollary~\ref{BUcontractible} we can conclude
\hbox {$|NU(x)| \simeq *$} as required.

\bigbreak

To prove~(2) it is enough to show that $\mathrm {hocolim}^\prime\, T(k)^x =
\emptyset$ for $x \in \bZ^n \cap \mathrm {int}\, kP$ and $\mathrm
{hocolim}^\prime\, T(k)^x \simeq *$ otherwise.  Since $\mathcal{O}_P(k)$ is
the same functor as~$\mathcal{O}_{-kP}(-1)$, it suffices to consider $k=-1$.

So assume $k = -1$.  Fix a point~$x \in \bZ^n$ and a face $F \in F(P)_0^1$.
Then $x \notin T(-1)^x (F)$ if and only if there is a facet $G \supseteq F$
of~$P$ such that $x$ and $\mathrm{int}\,(-P)$ are on the same side of the
affine hyperplane spanned by~$-G$. Such a facet certainly exists if $x \in
\mathrm{int}\,(-P)$. Hence $\mathrm {hocolim}^\prime\, T(k)^x = \emptyset$ in
this case.

If~$x$ is not in the interior of~$-P$, Lemma~\ref{LemmaFrontBack}, applied to
the polytope~$-P$, shows that $x \in T(-1)^x (F)$ if and only if $-F$~is a
front face of~$-P$ in the sense of Definition~\ref{DefFrontBack}. It follows
from Corollary~\ref{BFcontractible} that $\mathrm {hocolim}^\prime\, T(k)^x
\iso |NFront(x)| \simeq *$.
\end{proof}

\subsubsection*{Appendix: Cohomology of $X_P$}

The techniques from the computation of the space $\Gamma (\mathcal{O}_P (k))$
are applicable in the context of algebraic geometry: They can be used to give a
complete calculation of the cohomology groups of~$X_P$ with coefficients in a
twist of the structure sheaf. Let $R$ be a
commutative ring, and let $\mathcal{F}(k)$ denote the twisting sheaf $F
\mapsto \tilde R [\mathcal{O}_P (k)] = R [(C_F + kF) \cap \bZ^n]$. After
choosing orientations for the faces of~$P$, we can define a cochain complex
$C^\bullet$ of $R$-modules by setting $C^j := \bigoplus_{\dim F = j}
\mathcal{F}(k)^F$. The coboundary map is induced by
$$ \mathcal{F}(k)^F = R [(C_F + kF) \cap \bZ^n] \rTo[l>=4em]^{[F:G]}
R [(C_G + kG) \cap \bZ^n] = \mathcal{F}(k)^G $$
(for faces $F$,$G$ of~$P$ with $\dim G = 1+\dim F$) where $[F:G]$ is the incidence
number of $F$ and $G$. The cohomology groups of $C^\bullet$ are the cohomology
groups of~$X_P$ with coefficients in $\mathcal{F}(k)$
\cite[\S2]{H-Finiteness}: $H^r (X_P; \mathcal{F}(k)) \iso h^r (C^\bullet)$.

Now all the terms in $C^\bullet$ carry a natural $\bZ^n$-grading, and the
coboundary maps are homogeneous of degree~$0$.  Consequently, $C^\bullet$
splits into a direct sum of chain complexes $C^\bullet = \bigoplus_{x \in
  \bZ^n} C^\bullet_x$, and $h^k (C^\bullet) = \bigoplus_{x \in \bZ^n} h^k
(C^\bullet_x)$. The cochain complex $C^\bullet_x$ is given by
$$ C^j_x = \bigoplus_{{\scriptstyle \dim F = j} \atop
{\lower 2 pt \hbox{\({\scriptstyle x \in C_F + kF}\)}}} R $$
with coboundary maps given by incidence numbers as before.

Let $D^\bullet$ be the cochain complex defined by $D^j = \bigoplus_{\dim F =
j} R$ with coboundaries given by incidence numbers. Then $D^\bullet$ is the
cochain complex computing the (cellular) cohomology $H^* (P; R)$ of the
polytope~$P$. Hence $h^0 (D^\bullet) = R$, and $h^j (D^\bullet) = 0$ for $j
\not= 0$. Note that there is an inclusion map $C^\bullet_x \rTo D^\bullet$.

Now consider the case $k < 0$. If $x \in \mathrm{int}\, (kP)$, then the proof
of Theorem~\ref{HocolimT} shows that $x \notin C_F + kF$ for all proper faces
$F$ of~$P$, so $C^n_x = R$ and $C^j_x = 0 $ for $j \not= n$. Consequently, the
only non-vanishing cohomology group of $C^\bullet_x$ is $h^n(C^\bullet_x) =
R$.---If however $x \notin \mathrm{int}\, (kP)$, then the proof of
Theorem~\ref{HocolimT} shows that $x \in C_F + kF$ if and only if either $F=P$
or $F$ is a front face of~$P$ with respect to~$x$. We can thus identify the
quotient cochain complex $D^\bullet / C^\bullet_x$ with the cochain complex
computing the (cellular) cohomology of the space $|\Back(x)| \iso B^{n-1}$,
cf.~Lemma~\ref{FandBareBalls}. By the long exact sequence of cohomology groups
associated to the short exact sequence of cochain complexes
$$0 \rTo C^\bullet_x \rTo D^\bullet \rTo D^\bullet / C^\bullet_x \rTo 0$$
we infer that all cohomology groups of $C^\bullet_x$ are trivial.

Similar arguments apply to the cases $k=0$ and $k > 0$. By summation over all
$x \in \bZ^n$, we obtain:

\begin{theorem}
  Let $P \subseteq \bR^n$ be a lattice polytope with non-empty interior, let
  $R$ be a commutative ring, and denote by $\mathcal{F}(k)$ the quasi-coherent
  sheaf on~$X_P$ determined by $kP$. Let $k \in \bZ$ and $r \in \bN$.
  \begin{enumerate}
  \item If $k \geq 0$, then $H^r (X_P; \mathcal{F}(k)) = 0$ for $r \not= 0$, and
  there is an isomorphism $H^0 (X_P; \mathcal{F}(k)) \iso R [kP \cap
  \bZ^n]$. In particular $H^0 (X_P; \mathcal{F}(0)) \iso R$.
  \item If $k < 0$, then $H^r (X_P; \mathcal{F}(k)) = 0$ for $r \not= n$, and
  there is an isomorphism $H^n (X_P; \mathcal{F}(k)) = R [\mathrm{int}\, (kP)
  \cap \bZ^n]$. In particular $H^n (X_P; \mathcal{F}(k)) = 0$ if $E_P (k) = 0$.
  \end{enumerate}
  Thus the total cohomology $H^* (X_P; \mathcal{F}(k))$ is a free
  $R$-module of rank $|E_P (k)|$.\doqed
\end{theorem}

\begin{remark}
  \label{rem:compute_cohomology}
  As mentioned in the introduction of the paper, the interesting feature of
  this calculation is that it avoids the use of \textsc{Serre} duality
  \cite[\S7.7]{Danilov-toric} in favour of a topological argument. The reader
  might be interested in having a reference for the algebro-geometric version
  as well.  For $R$ a field, it follows from classical results in toric
  geometry \cite[Corollary~7.3]{Danilov-toric} that $H^r (X_P; \mathcal{F}(k))
  = 0$ for $r>0$ and $k \geq 0$, and that $H^0 (X_P; \mathcal{F}(k))$ has a
  canonical vector space basis given by the set $P \cap \bZ^n$,
  cf.~\cite[11.12]{Danilov-toric}. \textsc{Serre} duality implies that for
  negative~$k$ we have $H^r (X_P; \mathcal{F}(k)) = 0$ if $r\not= n = \dim
  (P)$.  Replacing $X_P$ by a non-singular variety and invoking \textsc{Serre}
  duality again, the argument given in \cite[\S11.12.4]{Danilov-toric}
  provides an alternative proof of the above theorem.
\end{remark}

\subsection{Computing $\psi_0 \circ \Gamma$}

We have calculated the composition $\Gamma \circ \psi_k$ in
Theorem~\ref{HocolimT} above. For the splitting result in $K$-theory we also
need to examine the composition $\psi_0 \circ \Gamma$: It is connected by a
chain of natural transformations to the functor $\Sigma^n$. We begin by
constructing two models $\epsilon$ and~$\sigma$ for the suspension functor on
the category of non-linear sheaves; the functor $\epsilon$ is naturally
isomorphic to $\Sigma^n$, and $\Sigma\sigma$ is naturally weakly equivalent
to~$\Sigma\epsilon$.

\begin{construction}
  Fix $Y \in \hcoh$. For a non-empty face~$A$ of~$P$, let $\mathrm{spr}_A Y$
  denote the diagram
  $$ \mathrm{spr}_A Y \colon F(P)_0 \rTo \Top,\quad
  F \mapsto Y^{A \vee F} \ .$$
  (In the language of algebraic geometry, $\mathrm{spr}_A Y$ describes the
  sheaf $f_* (Y|_{U_A})$ where $f \colon U_A = \mathrm{Spec}\,\tilde R[S_A]
  \subseteq X_P$ is the inclusion of an open affine. We retain the
  notation spr from \cite[\S3.8]{H-Proj} where $\mathrm{spr}_A Y$ is
  called a ``spread sheaf''.)  Since $A$ is a face of $A
  \vee F$ there is an inclusion of monoids $S_A \subseteq S_{A \vee F}$.
  Consequently, all spaces in the diagram $\mathrm{spr}_A Y$ have an
  $S_A$-action, so $\Gamma (\mathrm{spr}_A Y)$ has an $S_A$-action.  This
  construction is natural in~$A$: If $B$ is a face of~$P$ containing~$A$, the
  structure maps of~$Y$ define a natural transformation $\mathrm{spr}_A Y \rTo
  \mathrm{spr}_B Y$ with $S_A$-equivariant components. Consequently,
  application of~$\Gamma$ yields a diagram~$\sigma Y$, defined as
  $$ \sigma Y \colon F(P)_0 \rTo \Top,\quad A \mapsto \Gamma (\mathrm{spr}_A
  Y) \ .$$
  For a space $K \in \Top$ we define the constant diagram
  $$ \mathrm{con} K \colon F(P)_0 \rTo \Top,\quad F \mapsto K \ .$$
  Given $A \in F(P)_0$, the structure maps of~$Y$ define a natural
  transformation $\mathrm{con} Y^A \rTo \mathrm{spr}_A
  Y$. By naturality in~$A$ we obtain the diagram
  $$ \epsilon Y \colon F(P)_0 \rTo \Top,\quad F \mapsto \Gamma (\mathrm{con}
  Y^F) $$
  and a map of diagrams $\epsilon Y \rTo \sigma Y$.
\end{construction}

\begin{lemma}
\label{lem:sigma-is-suspension}
  \begin{enumerate}
  \item The diagram $\epsilon Y$ is naturally isomorphic to $\Sigma^n Y$. In
    particular, $\epsilon Y \in \hcoh$.
  \item The components of the map $\epsilon Y \rTo \sigma Y$ have
    contractible homotopy cofibres.
    In particular, the diagram $\Sigma (\sigma Y)$ is weakly
    equivalent to $\Sigma^{n+1} Y$ and thus is a non-linear sheaf.
  \item The functor~$\sigma$ defines a functor $\Sigma \sigma \colon \hcoh
    \rTo \wchcoh$.
  \end{enumerate}
\end{lemma}

\vglue 1 \smallskipamount

\begin{proof}
  (1).~For any space $K \in \C$ there are natural isomorphisms $\Gamma
  \mathrm{con} K \iso \Gamma \mathrm{con} (S^0 \wedge K) \iso (\Gamma
  \mathrm{con} S^0) \wedge K$. By definition of the total cofibre, the space
  $\Gamma \mathrm{con} S^0$ is the homotopy cofibre of the map $NF(P)_0^1 \rTo
  NF(P)_0$. Now $NF(P)_0^1$ is the barycentric subdivision of~$\partial P \iso
  S^{n-1}$, and $NF(P)_0$ is the barycentric subdivision of $P \iso B^n$.
  Consequently, $\Gamma \mathrm{con} S^0 \iso S^n$, proving the claim.

  (2).~The $A$-component of the map $\epsilon Y \rTo \sigma Y$ is
  given by applying $\Gamma$ to the natural transformation $\nu \colon
  \mathrm{con} Y^A \rTo \mathrm{spr}_A Y$. We want to show that the
  homotopy cofibre of~$\Gamma(\nu)$ is contractible. Since $\Gamma$
  commutes with taking homotopy cofibres, it is enough to show that
  the componentwise homotopy cofibre $Z = \mathrm{hocofibre}\,(\nu)$
  of the natural transformation~$\nu$ has contractible total cofibre.
  The diagram $Z \colon F(P)_0 \rTo \C$ is given as follows:
  \[Z^F = \mathrm{hocofibre}\, \left( Y^A = (\mathrm{con}Y^A)^F \rTo
    (\mathrm{spr}_A Y)^F = Y^{A \vee F} \right)\] We consider $Z$ as a
  diagram defined on~$F(P)$ by setting $Z^\emptyset = *$, and want to
  show $\Gamma (Z) \simeq *$.

  We claim that all the structure maps $Z^F \rTo Z^{A \vee F}$ are
  homotopy equivalences. In fact, for $F = \emptyset$ the source is a
  single point and the target is the homotopy cofibre of the identity
  on~$Y^A$ which is contractible. If $F \not= \emptyset$, the
  definition of~$Z$ shows that the structure map is the identity of
  the homotopy cofibre of $Y^A \rTo Y^{A \vee F}$.---By
  Theorem~\ref{GammaVanishes} this finishes the proof.

  (3).~This follows from~(2) and the fact that all constructions
  involved in the definition of~$\sigma$, when applied to locally (or
  even weakly) cofibrant objects, produce weakly cofibrant objects.
  Note that in general $\sigma Y$ will not be locally cofibrant; this
  happens, for example, if $Y = \mathcal{O}_P$.
\end{proof}

\medbreak

Recall that the structure maps of~$Y$ define a natural transformation of
diagrams $Y \rTo \mathrm{spr}_A Y$. The construction $\mathrm{spr}_A Y$ is
natural in~$A$, and taking total cofibres gives a map $\Gamma (Y) \rTo \lim_{A
\in F(P)_0} \Gamma \left( \mathrm{spr}_A Y \right) $.  The space
$\Gamma (\mathrm{spr}_A Y)$ has an $S_A$-action, and by passage to the adjoint (forcing
equivariance), we obtain a natural transformation
$$ \tau \colon \psi_0 \circ \Gamma (Y) \rTo \sigma Y $$
(where~$\psi_0$ is the canonical sheaf functor of Definition~\ref{def:canonical_sheaves}).

\begin{proposition}
\label{prop:psi0Gamma-is-suspension}
  The map of spaces $ \Sigma \Gamma (\tau) \colon \Sigma \Gamma \circ \psi_0
  \circ \Gamma (Y) \rTo \Sigma \Gamma (\sigma Y) $ is a homotopy equivalence.
\end{proposition}

\begin{proof}
  For this proof, we
  consider diagrams defined on $F(P)$, with $Y^\emptyset = *$. Let
  $\mathrm{spr}_\emptyset (Y)$ denote the trivial diagram with value~$*$
  everywhere. Then $(\mathrm{spr}_F Y)^G = (\mathrm{spr}_G Y)^F$ for all $F, G
  \in F(P)$. It is also convenient to define the pointed monoid $S_\emptyset =
  *$. We can now rewrite the map $\Gamma(\tau)$ as follows:

  \begin{diagram}
    \Gamma \left( F \mapsto \Gamma (Y) \wedge S_F \right) &
      \rTo^{\Gamma(\tau)} &
      \Gamma \left( F \mapsto \Gamma (\mathrm{spr}_F Y) \right) \\
    \dTo<\iso && \dTo<= \\
    \Gamma \left( F \mapsto \Gamma (G \mapsto Y^G \wedge S_F) \right) &
      \rTo[l>=6em]^{\Gamma(\Gamma (f))} &
      \Gamma \left(
        F \mapsto \Gamma (G \mapsto (\mathrm{spr}_F Y)^G) \right) \\
    \dTo<\iso>{\textrm{(Lemma~\ref{lem:GammaGamma})}} &&
      \dTo<\iso>{\textrm{(Lemma~\ref{lem:GammaGamma})}} \\
    \Gamma \left( G \mapsto \Gamma (F \mapsto Y^G \wedge S_F) \right) &
      \rTo^{\Gamma(\Gamma (g))} &
      \Gamma \left(
        G \mapsto \Gamma (F \mapsto (\mathrm{spr}_G Y)^F) \right) \\
  \end{diagram}
  Here the map $f$ is induced by the composition
  $$ Y^G \wedge S_F \rTo Y^{F \vee G} \wedge S_F
  \rTo[l>=3em]^{\mathrm{action}} Y^{F \vee G} =
  (\mathrm{spr}_F Y)^G \ , $$
  and~$g$ is the map $\psi_0 (Y^G) \rTo \mathrm{spr}_G Y$ with
  $F$-component (for the non-trivial case $G, F \not= \emptyset$)
  $$ g^F \colon (\psi_0 (Y^G))^F = Y^G \wedge S_F \rTo Y^{G \vee F} \wedge S_F
  \rTo[l>=3em]^{\mathrm{action}} Y^{G \vee F} = (\mathrm{spr}_G Y)^F \ .$$

  We will show that $\Gamma (g)$ suspends to a weak equivalence;
  then the same is true for $\Gamma (\Gamma (g)) \iso \Gamma (\tau)$
  as suspension commutes with taking total cofibres.

  We define two
  diagrams $R,Q \colon F(P) \rTo \Top$, depending on $G \in F(P)$, by setting
  \[R^F:= \cases { Y^G \wedge S^0   & if \( F = \emptyset \) \cr
                    (\psi_0 (Y^G))^F & if \( F \not= \emptyset \) \cr }
     \qquad \hbox{and} \qquad
     Q^F:= \cases { *   & if \( G = \emptyset \) \cr
                    Y^{G \vee F} & if \( G \not= \emptyset \) \cr }\]
  with structure map of~$R$ given by mapping the non-basepoint element
  of~$S^0$ to $0 \in S_F$ (the neutral element, not the base point). Structure
  maps in~$Q$ are the obvious ones induced by the structure maps
  of~$Y$. Except for their value at~$\emptyset$ the diagrams $R$ and~$Q$
  coincide with source and target of~$g$, respectively.

  Observe that $\Gamma (R)$ and $\Gamma (Q)$ are both contractible.  For the
  latter space this follows from Theorem~\ref{GammaVanishes} since all
  structure maps in $G$-direction are identities. For the former it follows
  from a slight modification of the proof of Theorem~\ref{HocolimT}. Note that
  $R = Y^G \wedge \overline\mathcal{O}_P (0)$ where $\overline\mathcal{O}_P (0)^F =
  \mathcal{O}_P (0)^F$ for all non-empty faces~$F$, and $\overline\mathcal{O}_P
  (0)^\emptyset = S^0=\{0\}_+$, the initial pointed monoid.  Thus it suffices
  to show $\Gamma (\overline\mathcal{O}_P (0)) \simeq *$.  In the notation used in
  Theorem~\ref{HocolimT} (the case $k=0$ and $x=0$), this means considering
  $T(0)^0$ as a functor on $F(P)^1$ with $T(0)^0 (\emptyset) = \{0\}$ whose
  (unpointed) homotopy colimit is $\mathrm{hocolim}^\prime\,T(0)^0 \iso
  |NF(P)^1| \simeq *$.

  For any space $K \in \Top$ let $\delta (K)$ denote the diagram which is
  trivial (with value~$*$) everywhere except that $\delta(K)^\emptyset = K$.
  We can build a commutative diagram
  \begin{diagram}
    \psi_0 (Y^G) & \rTo & R & \rTo & \delta (Y^G \wedge S^0) &
      \ \ \iso & \delta(Y^G)\\
    \dTo<g && \dTo && \dTo & \ldTo>{\delta(\id_{Y^G})} \\
    \mathrm{spr}_G Y & \rTo & Q & \rTo & \delta (Y^G) \\
  \end{diagram}
  where both rows are (componentwise) homotopy cofibre sequences. Indeed, the
  left horizontal maps are identities everywhere except possibly at
  $\emptyset$ in which case the source is a single point. So the natural maps
  from the homotopy cofibres to the diagrams on the right are weak
  equivalences.

  Applying $\Gamma$ to this diagram then gives a map of two cofibre sequences of
  topological spaces. By construction, the map on cofibres is the
  identity. The map on middle terms is a homotopy equivalence since $\Gamma(R)
  \simeq * \simeq \Gamma (Q)$ as remarked above. By considering the
  next step in the \textsc{Puppe} sequence of both rows we obtain a
  diagram of homotopy cofibre sequences
  \begin{diagram}
    \Gamma (R) & \rTo & \Gamma (\delta (Y^G)) & \rTo[l>=3em] & \Sigma \Gamma (\psi_0
    (Y^G)) \\
    \dTo<\simeq && \dTo<\id && \dTo>{\Sigma (\Gamma (g))}\\
    \Gamma (Q) & \rTo & \Gamma (\delta (Y^G)) & \rTo & \Sigma \Gamma (\mathrm{spr}_G
    Y) \\
  \end{diagram}
  which proves that $\Sigma \Gamma(g)$ is a weak equivalence.
\end{proof}

\section{Algebraic $K$-Theory of Non-Linear Sheaves}

For all of \S3, let $P \subset \bR^n$ be a lattice polytope with non-empty
interior.
\label{sec:algebraic-k-theory}

\subsection{Finiteness conditions}

\begin{definition}
  Let $Y$ be a non-linear sheaf on~$X_P$
  (Definition~\ref{Def:non-lin-sheaf}).
  \begin{enumerate}
  \item The object $Y$ is called {\it locally finite\/} if $Y^F \in
  \Cf(S_F)$, cf.~\S\ref{subsec:EquivSpaces},
    for all $F \in F(P)_0$. The full subcategory of~$\hcoh$
    (Definition~\ref{def:sheaves-cofibrancy}) consisting of
    locally finite non-linear sheaves is denoted $\hcohf$.
  \item The object $Y$ is called {\it homotopy finite\/} if it can be
    connected by a chain of weak equivalences to a locally finite non-linear sheaf.
    The full
    subcategory of~$\hcoh$ (Definition~\ref{def:sheaves-cofibrancy}) consisting of homotopy finite, locally
    cofibrant non-linear sheaves is denoted $\hcohhf$.
    The full subcategory of~$\wchcoh$ consisting of homotopy finite, weakly
    cofibrant (\ref{def:sheaves-cofibrancy}) non-linear sheaves is denoted
    $\wchcohhf$.
  \end{enumerate}
\end{definition}

\begin{remark}
  If a non-linear sheaf~$Y$ on~$X_P$ is homotopy finite then
  necessarily $Y^F \in \Chf(S_F)$ for all $F \in F(P)_0$. This latter
  condition is sufficient as well; in short, one chooses spaces
  $\tilde Z^F \in \Cf (S_F)$ and weak equivalences $\tilde Z^f \rTo
  Y^F$ and constructs, by induction on $\dim F$, a weak equivalence $Z
  \rTo Y$ with $Z \in \hcohf$. The components of~$Z$ will be built
  from the spaces~$\tilde Z^F$ by iterated mapping cylinder
  constructions used to strictify homotopy commutative diagrams. For
  $P$ a simplex a detailed argument is given in
  \cite[Lemma~4.1.2]{H-Proj}, the general case is similar.
\end{remark}

The canonical sheaf functors~$\psi_k$ defined
in~\ref{def:canonical_sheaves} preserve finiteness and weak
equivalences. Hence they restrict to functors $ \psi_k \colon \Chf
\rTo \hcohhf $.

\begin{proposition}
  The total cofibre construction restricts to a functor
  \[\Sigma \Gamma \colon \wchcohhf \rTo \Chf \ .\]
\end{proposition}

\begin{proof}
  For locally finite non-linear sheaves this is Theorem~3.9
  of~\cite{H-Finiteness}. Since both suspension and~$\Gamma$ are weakly
  homotopy invariant, the general case follows.
\end{proof}

\subsection{Algebraic $K$-theory and reduced $K$-theory}
\label{subsec:reduced-k-theory}

To define algebraic $K$-theory we use \textsc{Waldhausen}'s
$\mathcal{S}_\bullet$-construction for categories with cofibrations and weak
equivalences \cite{W-spaces}. We will work with the category $\hcohhf$ of
homotopy finite non-linear sheaves. A map $f \colon Y \rTo Z$ of non-linear
sheaves is called a {\it cofibration\/} if all its components are cofibrations
of equivariant spaces. The map $f$ is called an {\it $h$-equivalence\/} if it
is a weak equivalence, \ie, if all its components are weak homotopy
equivalences of spaces.  With respect to these cofibrations and weak
equivalences, we define the algebraic $K$-theory of the non-linear projective
toric variety~$X_P$ to be the space
$$ \Knl(X_P) := \Omega |h\mathcal{S}_\bullet \hcohhf| \ .$$

The functor $\Sigma \Gamma \colon \hcohhf \rTo \Chf$ is exact and thus
induces a map of $K$-theory spaces
$ \Knl(X_P) \rTo \Omega |h\mathcal{S}_\bullet \Chf| = A(*)$.
Roughly speaking the functor~$\psi_0$ provides a section up to homotopy of this
map; consequently we can split off a copy of $A(*)$ from $\Knl(X_P)$.

On a technical level, we use \textsc{Waldhausen\/}'s fibration theorem
\cite[Theorem~1.6.4]{W-spaces}. We call a map $f \colon Y \rTo Z$ of
non-linear sheaves an {\it $h_{[0]}$-equivalence\/} if
$\Sigma^2\Gamma(f)$ is a weak homotopy equivalence of spaces. Note
that the double suspension of the total cofibre of a non-linear sheaf
is a simply connected cofibrant pointed space, so {\it $f$ is an
  $h_{[0]}$-equivalence if and only if $\Gamma(f)$ induces an
  isomorphism of singular homology groups}. It follows that the class
of $h_{[0]}$-equivalences satisfy \textsc{Waldhausen}'s extension
axiom \cite[\S1.2]{W-spaces}. Since moreover every $h$-equivalence is
an $h_{[0]}$-equivalence, we can apply the fibration theorem to obtain
a fibration sequence
$$ \Omega |h \mathcal{S}_\bullet \hcohhf^{[0]}| \rTo^\iota
\Knl(X_P) \rTo
\Omega |h_{[0]} \mathcal{S}_\bullet \hcohhf| \eqno {(\dagger)} $$
where $\hcohhf^{[0]}$ is the full subcategory of~$\hcohhf$ consisting of
those objects~$Y$ satisfying $\Sigma^2 \Gamma(Y) \simeq *$ (\ie, the map $Y \rTo *$ is
an $h_{[0]}$-equivalence), and the map~$\iota$ is induced by
inclusion.

We need a lemma first. A map $f \colon Y \rTo Z$ of non-linear sheaves is
called a {\it weak cofibration\/} if all its components are cofibrations of
underlying pointed topological spaces.

\begin{lemma}
\label{lem:weak-cof-K}
  The inclusion $\hcohhf \subseteq \wchcohhf$ induces a weak equivalence
  $$ \Omega |h_{[0]}\mathcal{S}_\bullet \hcohhf| \simeq
  \Omega |h_{[0]}\mathcal{S}_\bullet \wchcohhf| $$
  where both $K$-theory spaces are defined with respect to
  $h_{[0]}$-equivalences, and on the right we use weak cofibrations.
\end{lemma}

\begin{proof}
  The category of diagrams $F(P)_0 \rTo \Top$ which satisfy conditions~(1)
  and~(2) of Definition~\ref{Def:non-lin-sheaf} has a \textsc{Quillen} closed
  model structure with cofibrations and weak equivalences ($h$-equivalences)
  as used for the category~$\hcoh$. This is a straightforward generalisation
  of \cite[Proposition~3.4.4]{H-Proj}, and can be considered as a special
  case of a model structure for ``twisted'' diagrams~\cite[Theorem~3.3.5]{HR-Twisted}.
  Consequently, every map $Y \rTo Z$ of a locally cofibrant object to a weakly
  cofibrant object can be factored as a cofibration $Y \rTo W $ (making $W$
  locally cofibrant) followed by a weak equivalence $ W \rTo Z$ (making $W$
  homotopy finite). Since a weak equivalence is an $h_{[0]}$-equivalence, we
  can now apply the Approximation Theorem \cite[1.6.7]{W-spaces}.
\end{proof}

\medbreak

We are now in the position to identify the base of the fibration
sequence~$(\dagger)$ with $A(*)$.

\begin{lemma}
  \label{lem:reducing}
  The functor $\Sigma^2 \Gamma \colon \hcohhf \rTo \Chf$ induces a homotopy equivalence
  $\Omega |h_{[0]} \mathcal{S}_\bullet \hcohhf| \simeq A(*)$.
\end{lemma}

\begin{proof}
  By Theorem~\ref{HocolimT}~(1) the composite $\Sigma^2 \Gamma \circ \psi_0$
  is weakly equivalent to $\Sigma^{n+2}$, hence induces a self homotopy
  equivalence of $A(*)$. Consequently, the map induced by~$\Sigma^2 \Gamma$ is
  surjective on homotopy groups.

  We want to show that the composition $\psi_0 \circ \Sigma^2\Gamma \iso
  \Sigma^2 \psi_0 \circ \Gamma$ is weakly equivalent, with respect to
  $h_{[0]}$-equivalences, to $\Sigma^{n+2}$. By
  Lemma~\ref{lem:weak-cof-K} it is enough to
  show that this is the case if both functors are considered as endofunctors
  on $\wchcohhf$.  By Lemma~\ref{lem:sigma-is-suspension} and
  Proposition~\ref{prop:psi0Gamma-is-suspension}, the functors $\Sigma^2
  \psi_0 \circ \Gamma$ and $\Sigma^{n+2}$ are connected by a chain of
  $h_{[0]}$-equivalences
  $$ \psi_0 \circ \Sigma^2 \Gamma \ \iso \ \Sigma^2 \psi_0 \circ \Gamma \rTo
     \Sigma^2 \sigma \lTo \Sigma^2 \epsilon \ \iso \ \Sigma^{n+2} \ ,$$
  thus induce self homotopy equivalences on the
  $K$-theory space $\Omega |h_{[0]} \mathcal{S}_\bullet \wchcohhf|$. In
  particular, the map induced by~$\Sigma^2 \Gamma$ is injective on homotopy
  groups. (Note that the chain of $h_{[0]}$-equivalences involves the
  functor~$\sigma$ which takes values in~$\wchcohhf$; this is the reason why
  weakly cofibrant objects are needed for the argument.)
\end{proof}

\begin{definition}
  The fibre of the fibration sequence~$(\dagger)$ is called the {\it reduced
  $K$-theory\/} of~$X_P$, written $\tilde \Knl(X_P)$.
\end{definition}

Thus~$(\dagger)$ yields a fibration sequence $ \tilde \Knl(X_P) \rTo^\iota
\Knl(X_P) \rTo[l>=3em]^{\Sigma^2 \Gamma} A(*) $.  Since $\Sigma^2\Gamma \circ
\psi_0$ induces a self homotopy equivalence of~$A(*)$ by
Theorem~\ref{HocolimT}~(1), $\psi_0$ provides a section up to
homotopy of the fibration sequence and we obtain a homotopy equivalence
$$ \tilde \Knl(X_P) \times A(*) \rTo[l>=3em]^{\iota \vee \psi_0}
\Knl(X_P) \ .$$

\subsection{Splitting $\tilde \Knl(X_P)$}

If the polytope $P$ does not contain lattice points in its interior it is
possible to split off further copies of~$A(*)$ from~$\tilde K(X_P)$. As in
\S\ref{subsec:reduced-k-theory}, this is done by producing suitable fibration
sequences with sections.

\begin{definition}
  For $k \in \bZ$ let $[k] := \{0,\,1,\,\ldots, \,k\}$. A map $f \colon Y \rTo
  Z$ of non-linear sheaves is called an $h_{[k]}$-equivalence if for all $ j
  \in [k] $ the map $\Sigma^2\Gamma (\theta_j(f))$ is a weak homotopy
  equivalence. (Here $\theta_j$ denotes the twisting functor of
  Definition~\ref{def:twisting}.) We denote by $\hcohhf^{[k]}$ the full
  subcategory of~$\hcohhf$ consisting of non-linear sheaves~$Y$ for which the
  map $Y \rTo *$ is an $h_{[k]}$-equivalence. The category $\wchcohhf^{[k]}$ is
  defined similarly as a full subcategory of $\wchcohhf$.
\end{definition}

\begin{lemma}
\label{lem:psi-factors}
  Let $k \geq 1$, and suppose $E_P (-k) = 0$, \ie, suppose that $kP$ does not
  contain lattice points in its interior. Then $\psi_{-k}$ can be considered as
  a functor $\psi_{-k} \colon \Chf \rTo \hcohhf^{[k-1]}$.
\end{lemma}

\begin{proof}
  Let $K \in \Chf$. We have to show that for each $\ell \in [k-1]$ the space
  $\Sigma^2 \Gamma (\theta_\ell(\psi_{-k}(K))) \iso \Sigma^2 \Gamma
  (\psi_{\ell-k} (K))$ is contractible. Now from Theorem~\ref{thm:Ehrhart}~(2)
  it is clear that $E_P(-k)=0$ implies $E_P(\ell-k)=0$ since $-k \leq \ell-k
  \leq -1$. Thus the claim follows from Theorem~\ref{HocolimT}~(2).
\end{proof}

\begin{lemma}
\label{lem:weak-cof-K2}
  For any $k \geq 1$, the inclusion $\hcohhf \subseteq \wchcohhf$ induces a weak
  equivalence
  $$ \Omega |h_{[k]}\mathcal{S}_\bullet \hcohhf^{[k-1]}| \simeq
  \Omega |h_{[k]}\mathcal{S}_\bullet \wchcohhf^{[k-1]}| $$
  where both $K$-theory spaces are defined with respect to
  $h_{[k]}$-equivalences, and on the right we use weak cofibrations.
\end{lemma}

\begin{proof}
  This is similar to the proof of Lemma~\ref{lem:weak-cof-K}.
\end{proof}

\begin{lemma}
\label{lem:base-is-a}
  Let $k \geq 1$, and suppose $E_P (-k) = 0$. The functor
  $$\Sigma^2 \Gamma \circ \theta_k \colon \hcohhf^{[k-1]} \rTo \Chf$$
  induces a homotopy equivalence
  $\Omega |h_{[k]} \mathcal{S}_\bullet \hcohhf^{[k-1]}| \simeq A(*)$.
\end{lemma}

\begin{proof}
  By Lemma~\ref{lem:psi-factors}, the functor $\psi_{-k}$ induces a map
  backwards. By Theorem~\ref{HocolimT}~(1), the composition
  $(\Sigma^2 \Gamma \circ \theta_k) \circ \psi_{-k} \iso \Sigma^2\Gamma \circ
  \psi_0$ is weakly equivalent to $\Sigma^{n+2}$, hence induces a self
  homotopy equivalence of $A(*)$. Consequently, the map induced by~$\Sigma^2
  \Gamma \circ \theta_k$ is surjective on homotopy groups.

  As in the proof of Lemma~\ref{lem:reducing} we see that the
  composition $\psi_0 \circ \Sigma^2\Gamma \iso \Sigma^2 \psi_0 \circ
  \Gamma $ is connected to $\Sigma^{n+2}$, both considered as an
  endofunctors of $\wchcohhf$, by a chain of $h_{[0]}$-equivalences.
  Consequently, the conjugate
  \[\theta_{-k} \circ \left( \psi_0 \circ
    \Sigma^2\Gamma \right) \circ \theta_k \iso \psi_{-k} \circ (
  \Sigma^2 \Gamma \circ \theta_k)\]
  is connected to~$\Sigma^{n+2}$ by a chain of $h_{[k]}$-equivalences
  (recall that, by definition, any object $Y \in \wchcohhf^{[k-1]}$
  has the property that $\Sigma^2 \Gamma (Y(\ell)) \iso \Sigma^2
  \Gamma \circ \theta_\ell (Y) \simeq *$ for all $\ell \in [k-1]$).
  Since the inclusion $\hcohhf \subseteq \wchcohhf$ induces an
  equivalence on $K$-theory spaces (Lemma~\ref{lem:weak-cof-K2}), it
  follows that $\psi_{-k} \circ (\Sigma^2\Gamma \circ \theta_k)$
  induces a self homotopy equivalence on the $K$-theory space $\Omega
  |h_{[k]} \mathcal{S}_\bullet \hcohhf^{[k-1]}|$. In particular, the
  map induced by~$\Sigma^2 \Gamma \circ \theta_k$ is injective on
  homotopy groups.
\end{proof}

\begin{lemma}
  Let $k \geq 1$, and suppose $E_P (-k) = 0$. The functor $\psi_{-k}$ induces
  a homotopy equivalence $ \Omega |h \mathcal{S}_\bullet \hcohhf^{[k]}| \times
  A(*) \rTo[l>=4em]^{\iota \vee \psi_{-k}} \Omega |h \mathcal{S}_\bullet
  \hcohhf^{[k-1]}|$. Here $\iota$ denotes the inclusion $\hcohhf^{[k]} \rTo
  \hcohhf^{[k-1]}$.
\end{lemma}

\begin{proof}
  By the Fibration Theorem \cite[Theorem~1.6.4]{W-spaces} there is a fibration
  sequence
  $$ \Omega | h \mathcal{S}_\bullet \hcohhf^{[k]} | \rTo
  \Omega | h \mathcal{S}_\bullet \hcohhf^{[k-1]} | \rTo
  \Omega | h_{[k]} \mathcal{S}_\bullet \hcohhf^{[k-1]} | \ .$$
  By Lemma~\ref{lem:base-is-a} the base space is homotopy equivalent
  to~$A(*)$, so we have a fibration sequence
  $$ \Omega | h \mathcal{S}_\bullet \hcohhf^{[k]} | \rTo
  \Omega | h \mathcal{S}_\bullet \hcohhf^{[k-1]} |
  \rTo[l>=4em]^{\Sigma^2\Gamma \circ
  \theta_k} A(*) \ .$$
  The functor $\psi_{-k}$ induces a section up to homotopy by
  Theorem~\ref{HocolimT}~(1) since $\Gamma \circ \theta_k \circ \psi_{-k}
  \iso \Gamma \circ \psi_0$.
\end{proof}

Note that for $k=1$ the target of this homotopy equivalence is nothing but
$\tilde\Knl(X_P)$. From these homotopy equivalences for
$k=1,\,2,\,\ldots$ we obtain:

\begin{theorem}
  \label{thm:splitting}
  Suppose $P \subset \bR^n$ is a lattice polytope with non-empty interior.
  Let $k$ denote the number of integral roots of the \textsc{Ehrhart}
  polynomial $E_P(T)$, and let $\iota$ denote the inclusion functor
  $\hcohhf^{[k]} \rTo \hcohhf$. The functor
  \[\hcohhf^{[k]} \times \underbrace{\Chf \times \ldots \times \Chf}_{(k+1)
    \hbox { factors}} \rTo[l>=8em]^{\iota \vee \psi_0 \vee \psi_{-1} \vee
    \ldots \vee \psi_{-k}} \hcohhf\]
  induces a homotopy equivalence
  \[\Omega |h \mathcal{S}_\bullet \hcohhf^{[k]}| \times \underbrace{A(*)
    \times \ldots \times A(*)}_{(k+1) \hbox{ factors}} \rTo^\sim \Knl(X_P) \ .
  \doqed\]
\end{theorem}

In the case $k=0$ this reduces to the splitting $\tilde \Knl(X_P)
\times A(*) \simeq \Knl(X_P)$. Note that the obstruction for
splitting off another copy of~$A(*)$ is the failure of
Lemma~\ref{lem:psi-factors}: The functor $\psi_{-k-1}$ will not
factorise through the category $\hcohhf^{[k]}$. In geometric terms,
this happens since the polytope $(k+1)P$ contains lattice points in
its interior (but $kP$ does not). In the language of algebraic
geometry, this means $H^* (X_P; \mathcal{F}(-k-1)) \ne 0$ while $H^*
(X_P; \mathcal{F} (-k)) = 0$, using the notation of the Appendix
to~\S\ref{subsec:total-cofib-can}.

\goodbreak

\end{document}